\documentclass[journal]{IEEEtran}
\usepackage[utf8]{inputenc}
\usepackage{color}
\usepackage{graphicx}          
\usepackage{cite}                               
\usepackage{amsmath}
\usepackage{mathpazo}
\usepackage{amsfonts}
\usepackage{amssymb}
\usepackage{mathrsfs}
\usepackage{url}
\usepackage{float}
\usepackage{lmodern}
\usepackage{algorithm}
\usepackage{algpseudocode}
\usepackage{amsmath}
\usepackage{amsthm}
\usepackage{amssymb}
\usepackage{graphicx}
\usepackage[unicode=true,
 bookmarks=false,
 breaklinks=false,pdfborder={0 0 1},colorlinks=false]
 {hyperref}
\hypersetup{
 colorlinks,bookmarksopen,bookmarksnumbered,citecolor=red,urlcolor=red}

\makeatletter
 \let\oldforeign@language\foreign@language
 \DeclareRobustCommand{\foreign@language}[1]{%
   \lowercase{\oldforeign@language{#1}}}
\theoremstyle{plain}

\theoremstyle{remark}

\theoremstyle{definition}

\theoremstyle{plain}

 \let\oldforeign@language\foreign@language
 \DeclareRobustCommand{\foreign@language}[1]{%
   \lowercase{\oldforeign@language{#1}}}
\ifCLASSINFOpdf
\else
\fi
\providecommand{\definitionname}{Definition}
\providecommand{\lemmaname}{Lemma}
\providecommand{\remarkname}{Remark}
\providecommand{\theoremname}{Theorem}

\providecommand{\definitionname}{Definition}
\providecommand{\lemmaname}{Lemma}
\providecommand{\remarkname}{Remark}
\providecommand{\theoremname}{Theorem}

\providecommand{\definitionname}{Definition}
\providecommand{\lemmaname}{Lemma}
\providecommand{\remarkname}{Remark}
\providecommand{\theoremname}{Theorem}

\providecommand{\definitionname}{Definition}
\providecommand{\lemmaname}{Lemma}
\providecommand{\remarkname}{Remark}
\providecommand{\theoremname}{Theorem}

\providecommand{\definitionname}{Definition}
\providecommand{\lemmaname}{Lemma}
\providecommand{\remarkname}{Remark}
\providecommand{\theoremname}{Theorem}

\providecommand{\definitionname}{Definition}
\providecommand{\lemmaname}{Lemma}
\providecommand{\remarkname}{Remark}
\providecommand{\theoremname}{Theorem}

\providecommand{\definitionname}{Definition}
\providecommand{\lemmaname}{Lemma}
\providecommand{\remarkname}{Remark}
\providecommand{\theoremname}{Theorem}

\providecommand{\definitionname}{Definition}
\providecommand{\lemmaname}{Lemma}
\providecommand{\remarkname}{Remark}
\providecommand{\theoremname}{Theorem}

\providecommand{\definitionname}{Definition}
\providecommand{\lemmaname}{Lemma}
\providecommand{\remarkname}{Remark}
\providecommand{\theoremname}{Theorem}

\providecommand{\definitionname}{Definition}
\providecommand{\lemmaname}{Lemma}
\providecommand{\remarkname}{Remark}
\providecommand{\theoremname}{Theorem}

\providecommand{\definitionname}{Definition}
\providecommand{\lemmaname}{Lemma}
\providecommand{\remarkname}{Remark}
\providecommand{\theoremname}{Theorem}

\providecommand{\definitionname}{Definition}
\providecommand{\lemmaname}{Lemma}
\providecommand{\remarkname}{Remark}
\providecommand{\theoremname}{Theorem}

\providecommand{\definitionname}{Definition}
\providecommand{\lemmaname}{Lemma}
\providecommand{\remarkname}{Remark}
\providecommand{\theoremname}{Theorem}

\providecommand{\definitionname}{Definition}
\providecommand{\lemmaname}{Lemma}
\providecommand{\remarkname}{Remark}
\providecommand{\theoremname}{Theorem}

\providecommand{\definitionname}{Definition}
\providecommand{\lemmaname}{Lemma}
\providecommand{\remarkname}{Remark}
\providecommand{\theoremname}{Theorem}

\providecommand{\definitionname}{Definition}
\providecommand{\lemmaname}{Lemma}
\providecommand{\remarkname}{Remark}
\providecommand{\theoremname}{Theorem}

\providecommand{\definitionname}{Definition}
\providecommand{\lemmaname}{Lemma}
\providecommand{\remarkname}{Remark}
\providecommand{\theoremname}{Theorem}

\providecommand{\definitionname}{Definition}
\providecommand{\lemmaname}{Lemma}
\providecommand{\remarkname}{Remark}
\providecommand{\theoremname}{Theorem}

\providecommand{\definitionname}{Definition}
\providecommand{\lemmaname}{Lemma}
\providecommand{\remarkname}{Remark}
\providecommand{\theoremname}{Theorem}

\providecommand{\definitionname}{Definition}
\providecommand{\lemmaname}{Lemma}
\providecommand{\remarkname}{Remark}
\providecommand{\theoremname}{Theorem}

\providecommand{\definitionname}{Definition}
\providecommand{\lemmaname}{Lemma}
\providecommand{\remarkname}{Remark}
\providecommand{\theoremname}{Theorem}

\providecommand{\definitionname}{Definition}
\providecommand{\lemmaname}{Lemma}
\providecommand{\remarkname}{Remark}
\providecommand{\theoremname}{Theorem}

\providecommand{\definitionname}{Definition}
\providecommand{\lemmaname}{Lemma}
\providecommand{\remarkname}{Remark}
\providecommand{\theoremname}{Theorem}

\makeatother

\providecommand{\definitionname}{Definition}
\providecommand{\lemmaname}{Lemma}
\providecommand{\remarkname}{Remark}
\providecommand{\theoremname}{Theorem}
   \newtheorem{remark}{Remark}
   \newtheorem{proposition}{Proposition}
   \newtheorem{definition}{Definition}
   \newtheorem{assumption}{Assumption}
   \newtheorem{property}{Property}
   \newtheorem{lemma}{lemma}
   \newtheorem{theorem}{Theorem}

\begin{document}


\title{Neuro-adaptive distributed control with prescribed performance for the synchronization of unknown nonlinear networked systems}

\author{Sami~El-Ferik, Hashim.~A.~Hashim$^*$, and~Frank L. Lewis
\thanks{Sami~El-Ferik is with the Department
of Systems Engineering, King Fahd University of Petroleum and Minerals, Dhahran, 31261 e-mail: selferik@kfupm.edu.sa}
\thanks{$^*$Corresponding author, H. A. Hashim is with the Department of Electrical and Computer Engineering,
University of Western Ontario, Ontario, Canada e-mail: hmoham33@uwo.ca.}
\thanks{Frank L. Lewis is with the Research Institute, The University of Texas at Arlington, Texas 76118 email:lewis@uta.edu}
}

\maketitle
\begin{abstract}
This paper proposes a neuro-adaptive distributive cooperative tracking control with prescribed performance function (PPF) for highly nonlinear multi-agent systems. PPF allows error tracking from a predefined large set to be trapped into a predefined small set. The key idea is to transform the constrained system into unconstrained one through transformation of the output error. Agents' dynamics are assumed to be completely unknown, and the controller is developed for strongly connected structured network. The proposed controller allows all agents to follow the trajectory of the leader node, while satisfying necessary dynamic requirements. The proposed approach guarantees uniform ultimate boundedness of the transformed error and the adaptive neural network weights. Simulations include two examples to validate the robustness and smoothness of the proposed controller against highly nonlinear heterogeneous networked system with time varying uncertain parameters and external disturbances.
\end{abstract}

\begin{IEEEkeywords}
Prescribed performance, Transformed error, Multi-agents, Neuro-Adaptive, Distributed adaptive control, Consensus, Transient, Steady-state error.
\end{IEEEkeywords}


\IEEEpeerreviewmaketitle{}

\section{Introduction}
    The coordination of animals in social groups such as ants foraging, birds swarming, fish schooling and so forth have attracted the attention of researchers to develop distributed collaborative control of multi-agents in engineering applications especially robotic systems. Collaborative work between agents empower the group to execute hard and large tasks whcih the individual agent cannot simply achieve on its own. For instance, they are capable of simplifying mission by dividing it in several complementary tasks, exchanging vital information that may save the group, increasing productivity, and in certain cases, imitating biological behaviors. Combining collaborative framework with autonomous vehicles could be beneficial in many applications such as surveillance, inspection, space explorations and in many other areas where the agents are required to work in a distributed manner and communicate using a particular network configuration. Hence, agents, often called nodes, are customized to follow at least one or more leaders. The network formed by such nodes is called a communication graph, and the vertices of the graph comprise of the nodes while its edges represent the connections between nodes. A graph may be undirected, and in such a case, there is no difference between the two vertices associated with an edge. On the other hand, it may be directed from one vertex to another, which indicates the direction of information flow between each of the nodes and its respective neighbors.

    In the literature, \cite{fax_information_2004} is one of the pioneer work that addressed the consensus in multi-agent systems. Consensus for passive nonlinear systems has been studied in \cite{chopra_passivity-based_2006}, and node consensus of cooperative tracking problem has been investigated in many research work such as \cite{lewis_cooperative_2013}, \cite{olfati-saber_consensus_2007}, \cite{ren_distributed_2008} and \cite{zhang_adaptive_2012}. Cooperative tracking control  was studied for a single node in \cite{das_distributed_2010} and \cite{cao_distributed_2012} and in for high order dynamics in \cite{zhang_adaptive_2012}. Neural-Network-based robust adaptive control has been addressed in \cite{hou2009decentralized} to solve the consensus problem of multiagent systems connecting by an undirected graph type network topology. In  \cite{cheng2010neural}, The same problem has been addressed but in the case of a network having a directed graph communication topology. also unlike \cite{cheng2010neural}, the control in \cite{hou2009decentralized} cannot solve the leader-following problem in case the leader possesses a time-varying state trajectory. \cite{das_distributed_2010} and \cite{zhang_adaptive_2012} developed a neuro-adaptive distributed control for digraph-connected heterogeneous agents having unknown nonlinear dynamics. In \cite{das_distributed_2010}, the authors considered nodes with single integrator and later on, in \cite{zhang_adaptive_2012}, high order systems have been addressed. Although, ,most of the previous studies assumed that the input function of the node dynamics is known. However, cooperative tracking control problems of systems with unknown input function were studied in \cite{theodoridis_direct_2012} and \cite{el-ferik_neuro-adaptive_2014}. Neuro-adaptive fuzzy was proposed to approximate unknown nonlinear dynamics and input functions, of which the centers of output membership functions were determined based on off-line trials \cite{theodoridis_direct_2012}. One common ground shared by these studies is that they all considered both the unknown nonlinear dynamics and the input function to be linear in parameters (LIP); see \cite{el-ferik_neuro-adaptive_2014,lewis_neural_1998} for typical instances and \cite{liu2016pulse} for pulse-modulated intermittent control. Ultimate stability of the tracking error was the main concern in all these studies.

    Model uncertainties and external disturbances affect considerably the closed loop performance of the controlled agents. In particular, the error is proven to be ultimately bounded and converge within a residual set having a size that depends on some unknown but bounded terms representing the model uncertainties and the external disturbances. The prediction of the steady-state or transient error behavior becomes almost impossible to establish analytically \cite{bechlioulis_robust_2008} and \cite{bechlioulis2014robust9}. Prescribed performance simply means that the tracking error is confined within an arbitrarily small residual set, and the convergence error should be within a given range. In addition, the convergence rate has to be less than a prescribed value and they should be less than a prescribed constant. Prescribed performance with robust adaptive control was mainly developed to provide a smooth control input signal, soft tracking, and to address the problem of accurate computation of the upper bounds for systematic convergence. In summary, adaptive prescribed performance-based control ensures convergence of the error within a prescribed bound, an overshot less than a prescribed value, a uniform ultimate boundedness property for the transformed output error, and a smooth adaptive tracking.

    Developing cooperative adaptive control for multi-agent systems based on prescribed performance will yield many merits. Prescribed performance means secluding the system tracking output error to start within large set and convert systematically into an arbitrarily small set with predefined range \cite{bechlioulis_robust_2008} and \cite{bechlioulis2014robust9}. The error with prescribed performance should have a convergence rate less than a predefined value and the maximum overshoot or undershoot should not exceed a certain limit. Prescribed performance with cooperative adaptive control will have the ability to increase the robustness of the control and reduce the control effort. The proper computations of upper and lower bounds of the prescribed performance functions provide soft tracking error with systematic convergence. Robust adaptive control with PPF was proposed in \cite{bechlioulis_robust_2008} to develop controllers for feedback linearizable systems. Since then, many papers considering prescribed performance used linearly parametrized neural network to approximate the unknown nonlinearities and disturbances (see for instance \cite{bechlioulis_adaptive_2009}, \cite{tong2015fuzzy} , \cite{li2015prescribed} , \cite{Li_2017}, \cite{karayiannidis2016model}, \cite{na_adaptive_2013} and \cite{yang_adaptive_2015}). The neuro-approximation-based adaptive control was developed for some particular system structures such as strict-feedback systems \cite{bechlioulis_adaptive_2009}, affine systems \cite{wang_verifiable_2010}, high order nonlinear systems \cite{kostarigka_adaptive_2012}. These papers with made various assumptions regarding the continuity of the input matrix. Using trial and error methods in defining neural weights, adaptive control with prescribed performance was proposed for MIMO uncertain chaotic systems using model reference adaptive control in \cite{mohamed_improved_2014}. Almost all of these studies considered a single autonomous systems.

    In this paper, we propose a decentralized neuro-adaptive cooperative control for a fleet of autonomous systems linked through a communication digraph. The different synchronization errors between all these nodes are expected to have dynamic responses that satisfy predefined characteristics determined by the designer. The dynamics of each node is subject to model uncertainties and to unknown but bounded external disturbances. The control scheme is based on prescribed performance for the transient as well as the steady-state dynamic performance for each node's synchronization error. To our knowledge, this has never been treated before. In this paper, the original prescribed performance scheme as proposed in \cite{bechlioulis_robust_2008} will be modified to take into considerations the interactions between nodes created through the consensus algorithm in the presence of random disturbances at each time instant. Hence, stable, non-oscillatory dynamics with bounded and smooth decentralized control inputs are guaranteed.
  \\ The rest of the paper is organized as follows. In Section \ref{Sec2}, preliminaries of graph theory and prescribed performance bounds will be presented. The problem formulation and local error synchronizations in the sense of PPF are presented in Section \ref{Sec3}. In Section \ref{Sec4}, we develop the control law formulation and we prove stability of the connected digraph. Section \ref{Sec5} includes simulation results which verify the robustness of the proposed control. Finally, we conclude and suggest future directions of research in Section \ref{Sec6}.\\
  \textbf{Notations:} The following notations are used throughout the paper.\\
    \begin{tabular}{lcl}
    $|\cdot|$ &:&
     absolute value of a real  number;\\
     $\|\cdot\|$ &:&
    Euclidean norm of a vector;\\
    $\|\cdot \|_F$ &:&  Frobenius norm
    of a matrix; \\
    ${ \rm Tr}\{\cdot\}$ &:&  trace of a matrix;\\
    $\sigma(\cdot)$ &:&  set of singular values of a matrix,  with
    the \\
    & & maximum singular value $\bar{\sigma}$ and \\
    & & the minimum singular value $\underline{\sigma}$;\\
    $P > 0 $ &:&
    indicates that the matrix $P$ is positive definite \\
    & & ($P \geq 0$ positive semi-definite;\\
    $\mathbb{I}_n$ &:&  identity matrix with $n$-dimension.\\
    ${\bf \underline{1}}$ &:&  unity vector $[1, \ldots, 1]^{\top} \ \in \mathbb{R}^n $
    where $n$ is the \\
    & & required appropriate dimension.
    \end{tabular}
\section{Preliminaries} \label{Sec2}
 \subsection{Basic graph theory}
 A graph is denoted by $\mathcal{G} = (\mathcal{V}, \mathcal{E})$ with a nonempty finite set of nodes (or vertices) $\mathcal{V} = \{\mathcal{V}_1, \mathcal{V}_2, \ldots, \mathcal{V}_n\}$, and a set of edges (or arcs) $E\subseteq \mathcal{V}\times \mathcal{V}$. $(\mathcal{V}_i,\mathcal{V}_j) \in E$ if there is an edge from node $i$ to node $j$. The topology of a weighted graph is often described by the adjacency matrix $A=[a_{ij}]\in \mathbb{R}^{N\times N}$ with weights $a_{ij} > 0$ if $(\mathcal{V}_{j}, \mathcal{V}_{i}) \in E$: otherwise $a_{ij} = 0$. Throughout the paper, the topology is fixed,
 i.e. $A$ is time-invariant, and the self-connectivity element $a_{ii} = 0$. A graph can be directed or undirected. A directed graph is called diagraph. The weighted in-degree of a node $i$ is defined as the sum of i-th row of $A$, i.e., $d_i=\sum_{j=1}^Na_{ij}$. Define the diagonal in-degree matrix $D ={\rm diag}(d_1,\ldots,d_N)\in \mathbb{R}^{N\times N}$ and the graph Laplacian matrix $L = D-A$. The set of neighbors of a node $i$ is $N_i = \{j|(\mathcal{V}_j\times \mathcal{V}_i)\in E\}$. If node $j$ is a neighbor of node $i$, then node $i$ can get information from node $j$ , but not necessarily vice versa. For undirected graph, neighborhood is a mutual relation. A direct path from node $i$ to node $j$ is a sequence of successive edges in the form $\{(\mathcal{V}_i,\mathcal{V}_k),(\mathcal{V}_k,\mathcal{V}_l), \ldots, (\mathcal{V}_m,\mathcal{V}_j\}$. A diagraph has a spanning tree, if there is a node (called the root), such that there is a directed path from the root to every other node in the graph. A diagraph is strongly connected, if for any ordered pair of nodes
 $[\mathcal{V}_i,\mathcal{V}_j]$ with $i\neq j$, there is  directed path from node $i$ to node $j$ \cite{ren_distributed_2008}).

\section{Problem Formulation} \label{Sec3}
    Consider the following nonlinear dynamics for the $i$th node
       \begin{equation}
         \label{eq:eq1}
         \begin{aligned}
           & \dot{x}_i = f_i\left(x_i\right) + u_i + w_i
         \end{aligned}
       \end{equation}
     where the state node is $x_i \in \mathbb{R}$, the control signal node $u_i \in \mathbb{R}$ and the unknown disturbance for each node is $w_i \in \mathbb{R}$. $f_i\left(x_i\right) \in \mathbb{R}$ is the unknown nonlinear dynamics and assumed to be Lipschitz. If we assumed agent state is a vector where $x_i \in \mathbb{R}^{n}, n > 1$ then each of $u_i$, $f_i\left(x_i\right)$ and $w_i$ $\in \mathbb{R}^{n}$ and further results will include Kronecker products. From \eqref{eq:eq1}, the global graph dynamics can be described by
       \begin{equation}
         \label{eq:eq2}
         \begin{aligned}
           & \dot{x} = f\left(x\right) + u + w
         \end{aligned}
       \end{equation}
     where $x = [x_1,\ldots,x_N]^{\top}\in \mathbb{R}^{N}$, $u = [u_1,\ldots,u_N]^{\top}\in \mathbb{R}^{N}$, $w = [w_1,\ldots,w_N]^{\top}\in \mathbb{R}^{N}$ and $f\left(x\right) = [f_1\left(x_1\right),\ldots,f_N\left(x_N\right)]^{\top}\in \mathbb{R}^{N}$. The leader state $x_0$ is defined by the desired synchronization trajectory and the leader nonlinear dynamics is described as
       \begin{equation}
         \label{eq:eq3}
         \begin{aligned}
           & \dot{x}_0 = f_0\left(x_0,t\right)
         \end{aligned}
       \end{equation}
     where $x_0 \in \mathbb{R}$ is leader node state and $f_0\left(x_0,t\right) \in \mathbb{R}$ is the nonlinear function of the leader.\\
     Defining the local neighborhood synchronization error function for node $i$ as in \cite{li_pinning_2004,khoo_robust_2009}
       \begin{equation}
         \label{eq:eq4}
         \begin{aligned}
           & e_i = \sum_{j \in N_i} a_{ij}(x_i-x_j)+b_{i}(x_i-x_0)
         \end{aligned}
       \end{equation}
     with pinning gains $a_{ij}\geq 0$ and $a_{ij} > 0$ for agent $i$ directed to the state of agent $j$, $b_{i}\geq 0$ and $b_{i}> 0$ for at least one agent $i$ is directed toward the leader state. From \eqref{eq:eq4}, the global error dynamics is described by
       \begin{equation}
         \label{eq:eq5}
         \begin{aligned}
          e &  = -\left(L+B\right)(\underline{x}_0-x) = \left(L+B\right)(x - \underline{x}_0)\\
            &  = \left(L+B\right)\;e_0
         \end{aligned}
       \end{equation}
     with $e = [e_1,\ldots,e_N]^{\top}\in \mathbb{R}^{N}$, $\underline{x}_0={\bf \underline{1}}x_0 \,\in\mathbb{R}^{N}$, $L \in \mathbb{R}^{N \times N}$, $B \in \mathbb{R}^{N \times N}$, $B ={\rm diag}\{b_i\}$ and ${\bf \underline{1}}=[1,\ldots,1]^{\top}\in\mathbb{R}^{N}$. Notice that $e_0 =x_0-x$. The proof of equation \eqref{eq:eq5} can be found in \cite{lewis_cooperative_2013}.\\
      \begin{remark}
      The communication graph is considered strongly connected. Thus, if $b_i \neq 0$ for at least one $i$, $i=1,\ldots, N$ then $\left(L+B\right)$ is an irreducible diagonally dominant M-matrix and hence nonsingular \cite{Qu2009}.
      \end{remark}

     \begin{remark} (see \cite{lewis_cooperative_2013})
         If agent state is $x_i \in \mathbb{R}^{n}$ and the leader state $x_0 \in \mathbb{R}^{n}$ where $n > 0$, then $e,x \in \mathbb{R}^{nN}$ and equation \eqref{eq:eq5} will be
     \end{remark}
       \begin{equation}
         \label{eq:eq6}
         \begin{aligned}
           & e = \left(\left(L+B\right)\otimes \mathbb{I}_N\right)(x - {\bf \underline{1}}x_0)
         \end{aligned}
       \end{equation}
     where $\otimes$ is the Kronecker product.
     The derivative error dynamics of \eqref{eq:eq5} is
       \begin{equation}
         \label{eq:eq7}
         \begin{aligned}
           & \dot{e} = \left(L+B\right)(f\left(x\right)+u+w-\underline{f}\left(x_0,t\right))
         \end{aligned}
       \end{equation}
     For strong connected graph, $B \neq 0$ and we have
       \begin{equation}
         \label{eq:eq8}
         \begin{aligned}
           & ||e_0|| \leq \frac{||e||}{\underline{\sigma}\left(L+B\right)}
         \end{aligned}
       \end{equation}
       where $\underline{\sigma}\left(L+B\right)$ is the minimum singular value of $L+B$.

In the case of multi-agent systems, Equation (\ref{eq:eq6}) reflects the coupling that has been created through the synchronization between the different states of each agent. Thus, these interactions no longer guarantee that the error dynamics of each agent will be confined within the desired performance functions, just based on knowledge of the sign $e_i(0)$, $i=1,...N$.

\subsection{Prescribed performance function} Prescribed performance function (PPF) has been defined as a time function tool of tracking error $e\left(t\right)$ such that $e\left(t\right)$ starts within a predefined large set and decays systematically to a predefined small set \cite{bechlioulis_robust_2008}, \cite{bechlioulis2014robust9} and \cite{mohamed_improved_2014}. PPF has been developed to provide smooth tracking response with sufficient range of control signal and it guarantees the transient and tracking performance with prescribed characteristics.\\
       A performance function $\rho\left(t\right)$ is associated with the error component $e\left(t\right)$ and is defined as a smooth function such that $\rho\left(t\right): \mathbb{R}_{+} \to \mathbb{R}_{+}$ is a positive decreasing function $\lim\limits_{t \to \infty}\rho\left(t\right)=\rho_{\infty}>0$. The PPF can be written as
       \begin{equation}
         \label{eq:eq9}
         \rho_i\left(t\right)=\left(\rho_{i0} - \rho_{i\infty}\right)\exp\left(-\ell_i\, t\right)+\rho_{i\infty}
       \end{equation}
       where $\rho_{i0},\rho_{i\infty}$ and $\ell_i$ are appropriately defined positive constants. In order to overcome the difficulty caused through the synchronization algorithm and achieve the desired prescribed performance, the following time varying constraints are proposed:
      \begin{equation}
         \label{eq:eq10b}
         -\delta_i\rho_i\left(t\right)<e_i\left(t\right)<\rho_i\left(t\right),\hspace{10pt} {\rm if} \: e_i\left(t\right)>0
       \end{equation}
       \begin{equation}
         \label{eq:eq11b}
         -\rho_i\left(t\right)<e_i\left(t\right)<\delta_i\rho_i\left(t\right),\hspace{10pt} {\rm if} \: e_i\left(t\right)<0
       \end{equation}
for all $ t \geq 0 $ and $ 0 \leq \delta_i \leq 1 $, and $i=1, ...,N$.

\begin{remark}
The dynamic constraints \eqref{eq:eq10b} and \eqref{eq:eq11b} represent a modification of the ones in \cite{bechlioulis_robust_2008},  \cite{wang2017dynamic} and \cite{mohamed_improved_2014}. In these papers, the constraints are conditioned on $e(0)$ as follows
\begin{equation}
         \label{eq:eq10}
         -\delta \rho \left(t\right)<e\left(t\right)<\rho \left(t\right),\hspace{10pt} {\rm if} \: e(0)>0
\end{equation}
       \begin{equation}
         \label{eq:eq11}
         -\rho \left(t\right)<e \left(t\right)<\delta \rho \left(t\right),\hspace{10pt} {\rm if} \: e(0)<0
       \end{equation}
 Due to the interaction between agents' dynamics such constraints will lead to instability. Figure \ref{fig:fig1} illustrates how the tracking error in the case of multi-agent systems may exceed the lower or upper bounds (plot in green color). Upon crossing this reference, the system becomes unstable under the original formulation \eqref{eq:eq10} and \eqref{eq:eq11}. However, the switching based on $e_i\left(t\right)$ provides the necessary control to keep the system stable. Figure~\ref{fig:fig1} shows the full idea of tracking error with prescribed performance as it transits from a large to a smaller set in accordance with equations \eqref{eq:eq10b} and \eqref{eq:eq11b}. Figure~\ref{fig:fig1} shows how the error is trapped systematically within large set to a predefined small set.
        \begin{figure}[h!]
         \centering
         \includegraphics[scale=0.27]{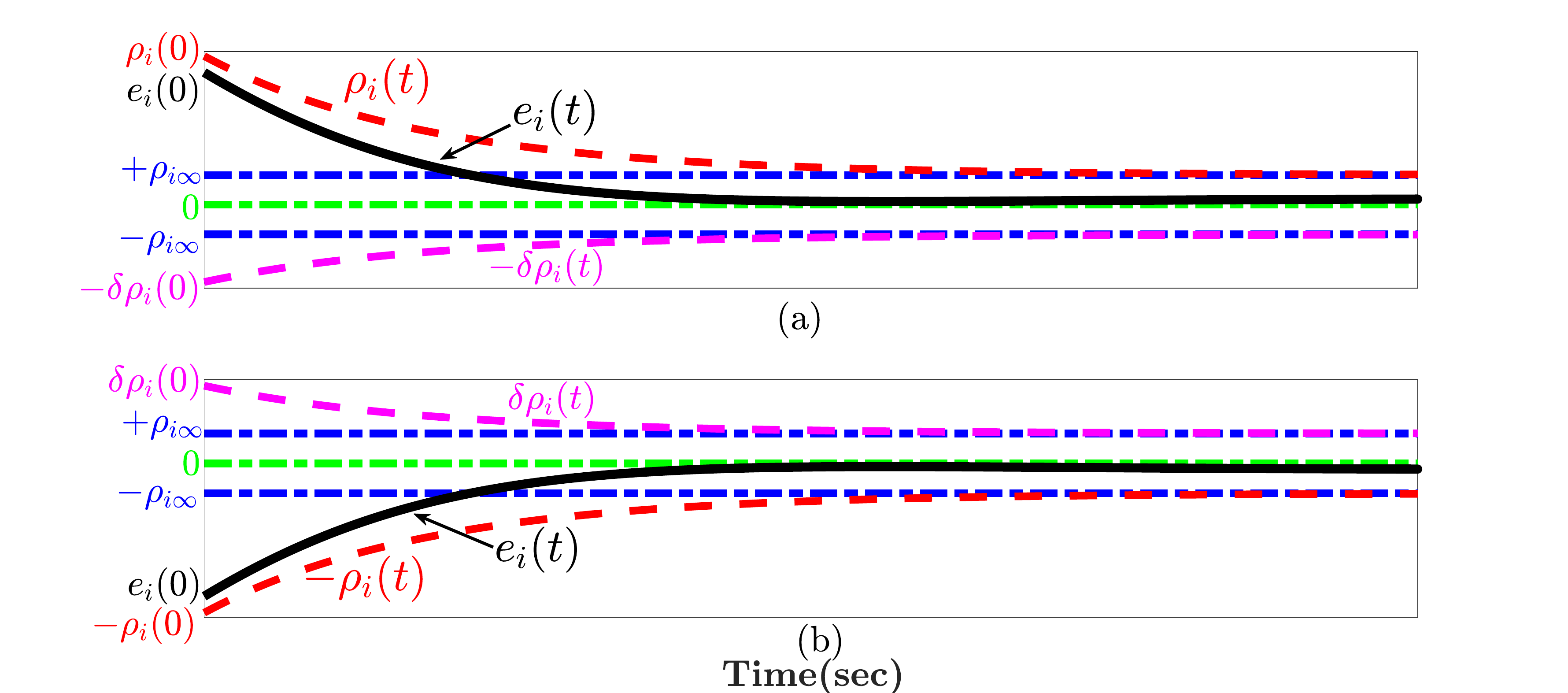}
         \caption{ Graphical representation of tracking error with prescribed performance
          (a) Prescribed performance of \eqref{eq:eq10b}; (b) Prescribed performance of \eqref{eq:eq11b}.}
         \label{fig:fig1}
        \end{figure}
\end{remark}
In order to transform the error of nonlinear system from the constraints in \eqref{eq:eq10b} and \eqref{eq:eq11b} to unconstrained, consider the following relation of transformed error given in \cite{bechlioulis_robust_2008}
      \begin{equation}
         \label{eq:eq12b}
         \epsilon_i=\psi\left(\frac{e_i\left(t\right)}{\rho_i\left(t\right)}\right)
      \end{equation}
      or equivalently,
      \begin{equation}
         \label{eq:eq13}
         e_i\left(t\right)=\rho_i\left(t\right)S\left(\epsilon_i\right)
      \end{equation}
      where $S(\cdot)$ is a smooth function and $\psi(\cdot)$ is its inverse for $i=1,2,\ldots,N$. In addition, $S(\cdot)$ is selected such that the following properties are satisfied:
      \begin{property}
      \label{propertiesS}
      \begin{enumerate}
        \item $S\left(\epsilon_i\right)$ is smooth and strictly increasing.
        \item $-\underline{\delta}_i<S\left(\epsilon_i\right)<\bar{\delta}_i,\hspace{10pt} {\rm if} \: e_i\left(t\right) \geq 0$\\
          $-\bar{\delta}_i<S\left(\epsilon_i\right)<\underline{\delta}_i,\hspace{10pt} {\rm if} \: e_i\left(t\right)<0$
        \item
            $\left.
             \begin{aligned}
               {\rm lim}_{\epsilon_i \rightarrow -\infty}S\left(\epsilon_i\right)=-\underline{\delta}_i\\
               {\rm lim}_{\epsilon_i \rightarrow +\infty}S\left(\epsilon_i\right)=\bar{\delta}_i
             \end{aligned}
             \right\}
             \qquad {\rm if} \hspace{10pt} e_i\left(t\right)\geq 0$\\
            $\left.
             \begin{aligned}
               {\rm lim}_{\epsilon_i \rightarrow -\infty}S\left(\epsilon_i\right)=-\bar{\delta}_i\\
               {\rm lim}_{\epsilon_i \rightarrow +\infty}S\left(\epsilon_i\right)=\underline{\delta}_i
             \end{aligned}
             \right\}
             \qquad {\rm if} \hspace{10pt} e_i\left(t\right)< 0$
         \end{enumerate}
     \end{property}
      where $\underline{\delta}_i$ and $\bar{\delta}_i$ are known positive constants that have to be selected according to the following relations:  \\
      \begin{equation} \label{eq:eq14b}
      \begin{split}
          & S\left(\epsilon_i\right)= \\
          & \left\{
          \begin{aligned}
          \frac{\bar{\delta}_i\exp\left(\epsilon_i\right)-\underline{\delta}_i\exp\left(-\epsilon_i\right)}{\exp\left(\epsilon_i\right)+\exp\left(-\epsilon_i\right)},& \quad {\rm with}\;\; \bar{\delta}_i> \underline{\delta}_i \hspace{5pt} {\rm if} \hspace{5pt} e_i\left(t\right)\geq 0\\
          \frac{\bar{\delta}_i\exp\left(\epsilon_i\right)-\underline{\delta}_i\exp\left(-\epsilon_i\right)}{\exp\left(\epsilon_i\right)+\exp\left(-\epsilon_i\right)},& \quad {\rm with} \;\; \underline{\delta}_i > \bar{\delta}_i \hspace{5pt} {\rm if} \hspace{5pt} e_i\left(t\right) < 0\\
          \end{aligned}
         \right.
         \end{split}
      \end{equation}
      Now, consider the general form of the smooth function
      \begin{equation}
        \label{eq:eq15b}
        S\left(\epsilon_i\right)=
          \begin{aligned}
         & \frac{\bar{\delta}_i\exp\left(\epsilon_i\right)-\underline{\delta}_i\exp\left(-\epsilon_i\right)}{\exp\left(\epsilon_i\right)+\exp\left(-\epsilon_i\right)}
          \end{aligned}
      \end{equation}
     and the transformed error
      \begin{equation}
        \label{eq:eq16b}
         \begin{aligned}
         \epsilon_i = & S^{-1}\left(\frac{e_i\left(t\right)}{\rho_i\left(t\right)}\right)\\
          = & \frac{1}{2}\left\{ \begin{aligned} & {\rm ln} \frac{\underline{\delta}_i+e_i\left(t\right)/\rho_i\left(t\right)}{\bar{\delta}_i-e_i\left(t\right)/\rho_i\left(t\right)}, &  {\rm with}\;\; \bar{\delta}_i > \underline{\delta}_i \hspace{5pt} {\rm if} \hspace{5pt} e_i\left(t\right)\geq 0 \\
         & {\rm ln} \frac{\underline{\delta}_i+e_i\left(t\right)/\rho_i\left(t\right)}{\bar{\delta}_i-e_i\left(t\right)/\rho_i\left(t\right)} ,&{\rm with}\;\; \underline{\delta}_i > \bar{\delta}_i \hspace{5pt} {\rm if} \hspace{5pt} e_i\left(t\right)< 0\\ \end{aligned}  \right.
         \end{aligned}
       \end{equation}
   One can notice in the previous set of equations, $\underline{\delta}_i$ and $\bar{\delta}_i$ exchange values depending on the sign of $e_i\left(t\right)$. One should note that the highest value of both involves subtracting the absolute value of $e_i\left(t\right)/\rho_i\left(t\right)$ and the lowest involves the addition of the absolute value of $e_i\left(t\right)/\rho_i\left(t\right)$. Accordingly and recalling that $\rho_i\left(t\right)>0$, equations (\ref{eq:eq16b}) $\bar{\delta}_i > \underline{\delta}_i$ can be rewritten as
      \begin{equation}
        \label{eq:eq56b}
         \begin{aligned}
         \epsilon_i =
          & \frac{1}{2}\left\{ \begin{aligned} & {\rm ln} \frac{\underline{\delta}_i+ |e_i\left(t\right)|/\rho_i\left(t\right)}{\bar{\delta}_i-|e_i\left(t\right)|/\rho_i\left(t\right)}, &  \hspace{1pt} {\rm if} \hspace{1pt} e_i\left(t\right)\geq 0 \\
         & {\rm - ln} \frac{\underline{\delta}_i+|e_i\left(t\right)|/\rho_i\left(t\right)}{\bar{\delta}_i-|e_i\left(t\right)|/\rho_i\left(t\right)} ,&  \hspace{1pt} {\rm if} \hspace{1pt} e_i\left(t\right)< 0\\ \end{aligned}  \right.
         \end{aligned}
       \end{equation}

     Thus, the transformed error in \eqref{eq:eq12b} can be expressed in more compact form with $\bar{\delta}_i > \underline{\delta}_i$ as follows:
      \begin{equation}
        \label{eq:eq57b}
         \begin{aligned}
         \epsilon_i =  \frac{1}{2}{\rm sign}\left(e_i\left(t\right)/\rho_i\left(t\right)\right)\cdot {\rm ln} \left(\frac{\underline{\delta}_i+ |e_i\left(t\right)|/\rho_i\left(t\right)}{\bar{\delta}_i-|e_i\left(t\right)|/\rho_i\left(t\right)}\right)
         \end{aligned}
       \end{equation}
     And to attenuate the effect of chattering, the following form of the transformed error is proposed
      \begin{equation}
        \label{eq:eq58b}
         \begin{aligned}
         \epsilon_i =  \frac{1}{2\sqrt{\pi}}erf\left(\frac{\xi e_i\left(t\right)}{\rho_i\left(t\right)}\right)\cdot {\rm ln} \left(\frac{\underline{\delta}_i+ |e_i\left(t\right)|/\rho_i\left(t\right)}{\bar{\delta}_i-|e_i\left(t\right)|/\rho_i\left(t\right)}\right)
         \end{aligned}
       \end{equation}
     where $erf(\xi e/\rho)=\frac{2}{\sqrt{\pi}}\int\limits_{0}^{\xi e}e^{-a^2}da$. $\xi>0$ is a design parameter.
     \begin{remark} The primary role of $\xi$ is to make $erf(\xi e)$ as close as possible to ${\rm sign}(e)$. Ideally $\xi$ is selected as big as possible. For instance, $|erf(\xi e)|\approxeq 1$  when $ |\frac{e}{\rho}| \geq\Delta=\frac{2}{\xi}$. Therefore, if $\xi=200$ then $|erf(e/\rho)|\approxeq 1$ when $|e/\rho|\geq0.01$.  However, while the error function derivative is smooth and bounded, the more one selects a high gain  $\xi$, the more the risk of chattering.
     \end{remark}
     For simplification, let's define $x:=x\left(t\right)$, $e:=e\left(t\right)$, $\epsilon:=\epsilon\left(t\right)$ and $\rho:=\rho\left(t\right)$. Then, after algebraic manipulations, the derivative of transformed error when $|e|/\rho \geq \Delta/\xi$ can be approximated by:
       \begin{equation}
        \label{eq:eq17}
       \dot{\epsilon}_i = \frac{1}{2\rho_i}\left(\frac{1}{\underline{\delta}_i+|e_i|/\rho_i} + \frac{1}{\bar{\delta}_i - |e_i|/\rho_i}\right)\left(\dot{e}_i - \frac{e_i\dot{\rho}_i}{\rho_i}\right)
       \end{equation}
 \begin{remark}
     As mentioned earlier, the selection of the high gain $\xi$ can make the absolute value of the error function converge to 1 for a very small ratio $|e\left(t\right)|/\rho=\Delta/\rho$. In our analysis, we will use \eqref{eq:eq17} to show that the control will generate a UUB error dynamic that will converge to a ball around zero with a radius that can be made as small as desired depending on the selection of $\xi$. Thus, the error may not converge to zero.
\end{remark}
     Let
      \begin{equation}
           \label{eq:eq17c}
            r_i = \frac{1}{2\rho_i}\left(\frac{1}{\underline{\delta}_i+|e_i|/\rho_i} + \frac{1}{\bar{\delta}_i - |e_i|/\rho_i}\right)
       \end{equation}
     From \eqref{eq:eq7} and \eqref{eq:eq17}, the global synchronization of the transformed error can be obtained as
       \begin{equation}
        \label{eq:eq18}
       \dot{\epsilon} = R\left(L+B\right)(f\left(x\right)+u+w-\underline{f}\left(x_0,t\right)) - \,\dot{\varUpsilon}\,\varUpsilon^{-1}\,e\left(t\right)
       \end{equation}
       where the control at the level of each node is of the form $u_i = -c \epsilon_i + \nu$; the value of $\nu$ will be clarified later in the paper, and it represents the necessary control actions to tackle the uncertainties and the Neural-Network approximation errors (see \ref{eq:eq27});
      $\epsilon=[\epsilon_1,\ldots,\epsilon_N]^{\top}\in \mathbb{R}^N $, $\varUpsilon={\rm diag}\left[\rho_i\left(t\right)\right]$ and $\dot{\varUpsilon}={\rm diag}\left[\dot{\rho_i}\left(t\right)\right]$, $i=1, \ldots,N$;
$R$ is such that $R = {\rm diag}\left[r_1\left(t\right),\ldots,r_N\left(t\right)\right]$ with  $R>0$ and $\dot{R}<0$;  $\,\dot{\varUpsilon}\,\varUpsilon^{-1}<0$ with ${\rm lim}_{t \to \infty}\,\dot{\varUpsilon}\,\varUpsilon^{-1} = 0$. It should be noted that $\bar{\delta}_i$ and $\underline{\delta}_i$  define the dynamic boundaries of the initial (large) and final (small) set and they have a significant impact on the control effort up to small error tracking  which is bounded by $\rho_{i\infty}$. Higher values of $\bar{\delta}_i$ and $\underline{\delta}_i$ require more time for the systematic convergence from large to small set. Whereas, the impact of $\ell_i$ and $\xi$ can be noticed on the speed of convergence from large to small sets that implies values of $\ell_i$ and $\xi$ have a direct impact on the range of control signal. Before proceeding further, the following definitions are needed (see \cite{das_distributed_2010}).
 \begin{definition}
 The global neighborhood error $e\left(t\right)\in \mathbb{R}^N $ is uniformly ultimately bounded (UUB) if there exists a compact set $\Omega \subset \mathbb{R}^N $ so that $\forall e\left(t_0\right) \in \Omega$ there exists a bound $B$ and a time $t_f\left(B,e\left(t_0\right)\right)$, both independent of $t_0 \geq 0$, such that $\left\Vert e\left(t\right) \right\Vert \leq B$ so that $\forall t > t_0+t_f$.
 \end{definition}
 \begin{definition}
 The control node trajectory $x_0\left(t\right)$ given by \eqref{eq:eq1} is cooperative UUB with respect to solutions of node dynamics \eqref{eq:eq3} if there exists a compact set $\Omega \subset \mathbb{R}^N $ so that $\forall \left(x_i\left(t_0\right)-x_0\left(t_0\right)\right) \in \Omega$, there exist a bound $B$ and a time $t_f\left(B,\left(x\left(t_0\right) - x_0\left(t_0\right)\right)\right)$, both independent of $t_0\geq 0$, such that $\left\Vert x\left(t_0\right) - x_0\left(t_0\right) \right\Vert \leq B$, $\forall i$, $\forall t > t_0+t_f$.
 \end{definition}
 \section{Neuro Adaptive Distributed Control with Prescribed Performance} \label{Sec4}
 This section presents the neural approximation of unknown nonlinearities and model uncertainties of a group of heterogeneous agents. The distributed control design is formulated for each local node to satisfy pre-defined prescribed characteristics.
 \subsection{Neural Approximations} 
     The unknown nonlinearities of local agents in \eqref{eq:eq1} can be approximated by
       \begin{equation}
        \label{eq:eq19}
           f_i\left(x_i\right) = W_i^{\top}\phi_i\left(x_i\right) + \alpha_i
       \end{equation}
     with $\phi_i\left(x_i\right) \in \mathbb{R}^{v_i}$ where $v_i$ is a sufficient number of neurons at each node, $W_i \in \mathbb{R}^{v_i}$ and $\alpha_i$ is the approximated error. According to \cite{hornik_multilayer_1989,lewis_neural_1998}, the neural network can be approximated by a variety of sets including radial basis functions with centers and widths \cite{poggio_regularization_1990}, sigmoid functions \cite{cotter_stone-weierstrass_1989}, etc.\\
     Tracking the local performance of each node will be attained through the compensation of unknown nonlinearities and using the available information of neighbor state agents. Therefore, the nonlinearities of local nodes can be approximated such as
       \begin{equation}
        \label{eq:eq20}
           \hat{f}_i\left(x_i\right) = \hat{W}_i^{\top}\phi_i\left(x_i\right)
       \end{equation}
        where $\hat{W}_i \in \mathbb{R}^{v_i}$ and $\hat{f}_i\left(x_i\right)$ is the approximation of $f_i\left(x_i\right)$. One should mention a crucial and well known property of the NN structure adopted in (\ref{eq:eq20})  (for more details the reader is invited to see \cite{lewis_neural_1998}).
               \begin{property}
               \label{property1}
               For each continuous function $f_i: \mathbb{R}^n\longrightarrow\mathbb{R}$ and for every $\epsilon_i>0$, there exist a bounded integer $v_i$, and  optimal synaptic weight vector $W_i^*\in \mathbb{R}^{v_i}$, $i=1,\ldots, N$ such that
               \begin{equation*}
               {\rm sup}_{x \in \Omega_x} \left| f\left(x\right) - W_i^{*T} \phi_i\left(x_i\right) \right|\leq \epsilon_i
               \end{equation*}
               \end{property}
     where $\Omega_x$ is a compact set. The global synchronization of $f\left(x\right)$ for the graph $\mathscr{G}$ can be described as
       \begin{equation}
        \label{eq:eq21}
           f\left(x\right) = W^{\top}\phi\left(x\right) + \alpha
       \end{equation}
       with $W^{\top}={\rm diag}\{W_i\},i=1,\ldots,N$, $\phi\left(x\right)=[\phi_1\left(x_1\right),\ldots,\phi_N\left(x_N\right)]$, $\alpha=[\alpha_1,\ldots,\alpha_N]$ and the global estimate of $f\left(x\right)$
       \begin{equation}
        \label{eq:eq22}
           \hat{f}\left(x\right) = \hat{W}^{\top}\phi\left(x\right)
       \end{equation}
       and $\hat{W}^{\top}={\rm diag}\{\hat{W}_i\},i=1,\ldots,N$. The error of estimating nonlinearities can be described by
        \begin{equation}
         \label{eq:eq23}
            \tilde{f}\left(x\right) = f\left(x\right) - \hat{f}\left(x\right) =\tilde{W}^{\top}\phi\left(x\right) + \alpha
        \end{equation}
        where $\tilde{W} = W - \hat{W}$ is the parameter estimation error.
   Owing to property \ref{property1}, there exist an unknown constant $\alpha_{Mi}>0$, such that $|\alpha_i| \leq \alpha_{Mi}$ for $i=1,\ldots,N$.

\subsection{NN Adaptive Control Design of Distributed Agents}
  As in \cite{das_distributed_2010}, the following standard assumptions are required.
         \begin{assumption}
                \label{assumption1}~~\\
                 a. Leader states are bounded by $\left\Vert x_0 \right\Vert \leq X_0$.\\
                 b. Leader unknown variable dynamics in \eqref{eq:eq3} are bounded such as $\left\Vert \underline{f}\left(x_0,t\right) \right\Vert \leq F_M$.\\
                 c. Unknown disturbances are bounded by $\left\Vert w \right\Vert \leq w_M$.\\
                 d. The variable $\phi$ in \eqref{eq:eq21} and \eqref{eq:eq22} is bounded by $\left\Vert \phi \right\Vert \leq \phi_M$.
          \end{assumption}
          \begin{lemma}
                \label{lemma1} 
                Let $L$ be an irreducible matrix and $B \neq 0 $ such as $\left(L+B\right)$ is nonsingular then we can define
                \begin{equation}
                 \label{eq:eq24}
                    q = [q_1,\ldots,q_N]^{\top}=\left(L+B\right)^{-1}\cdot{\bf \underline{1}}
                \end{equation}
                \begin{equation}
                 \label{eq:eq25}
                    P ={\rm diag}\{p_i\}={\rm diag}\{1/q_i\}
                \end{equation}
                Then, $P > 0$ and the matrix $Q$ defined as
                \begin{equation}
                 \label{eq:eq26}
                 \begin{aligned}
                    Q &= P\,\left(L+B\right)+\left(L+B\right)^{\top}\,P\\
                    &=P\left[\mathcal{K}\left(L+B\right)+\left(L+B\right)^{\top}\,\mathcal{K}\right]P
                                     \end{aligned}
                \end{equation}
                is also positive definite with  $\mathcal{K}:=P^{-1}$.
           \end{lemma}
           The gist of the idea is that $Q=\mathcal{K}\left(L+B\right)+\left(L+B\right)^{\top}\,\mathcal{K}$ is diagonally strictly dominant and since it is a symmetric M-matrix, then it is positive definite. Based on this lemma, the following Preposition holds.
           \begin{proposition}
            Let $R$ a positive definite diagonal matrix, and $L$, $B$, $P$ and $\mathcal{K}$ as defined in Lemma 1, then the matrix $Q$ defined as
                            \begin{equation}
                             \label{eq:eq261}
                                Q = P\,R\,\left(L+B\right)+\left(L+B\right)^{\top}\,R\,P
                            \end{equation}
            is positive definite.
           \end{proposition}
           \textbf{Proof:}\\
           Since $\left(L+B\right)$ is is a nonsingular M-matrix and $R>0$ is diagonal, then $R\left(L+B\right)$ is a non-singular M-Matrix.
           \begin{equation}
             \left(L+B\right)\,q ={\bf \underline{1}}>0
           \end{equation}
           Let $\mathcal{K}={\rm diag}\{q_i\}$ then
           \begin{equation}
           R\left(L+B\right)\mathcal{K} {\bf \underline{1}} = R\,\left(L+B\right)\,q=R{\bf \underline{1}}>0
           \end{equation}
          which means strict diagonal dominance of $R\,\left(L+B\right)\,\mathcal{K}$.
          \begin{equation}\label{eq:eq266}
          \begin{aligned}
            Q = & P\,R\,\left(L+B\right)+\left(L+B\right)^{\top}\,R\,P\\
              = & P\left[R\left(L+B\right)\mathcal{K}+\mathcal{K}\,\left(L+B\right)^{\top}\,R\right]P
          \end{aligned}
           \end{equation}
           $\left[R\left(L+B\right)\mathcal{K}+\mathcal{K}\,\left(L+B\right)^{\top}\,R\right]$ is symmetric and strictly diagonally dominant. Therefore,  $Q$ is positive definite.

          The control signal of local nodes is given by
                 \begin{equation}
                  \label{eq:eq27}
                     u_i = -c\epsilon_i - \hat{W_i}^{\top}\phi_i\left(x_i\right)
                 \end{equation}
                 with the control gain $c>0$. Let the control signal for all nodes be
                 \begin{equation}
                  \label{eq:eq28}
                     u = -c\epsilon - \hat{W}^{\top}\phi\left(x\right)
                 \end{equation}
                Let the NN local node tuning laws given by
                 \begin{equation}
                  \label{eq:eq29}
                     \dot{\hat{W}}_i = F_i\phi_i\epsilon_ir_ip_i(d_i+b_i) - kF_i\hat{W}_i
                 \end{equation}
                with $F_i \in \mathbb{R}^{v_i \times v_i}$, $F_i := \Pi_i\mathbb{I}_{v_i}$, $\Pi_i >0$ and $k \,>\,0$ are scalar gains.
                
           \begin{theorem}{\bf Distributed Adaptive Control Protocol for Synchronization with Prescribed Performance.}\\
           Consider the strong connected graph of the networked system in \eqref{eq:eq1} under Assumption \ref{assumption1} with distributed adaptive control protocol \eqref{eq:eq27} and  neighborhood synchronization errors given in \eqref{eq:eq4}. Let the adaptive estimate be defined as in \eqref{eq:eq29} and the local node NN tuning laws follow  $F_i = \Pi_i\mathbb{I}_{v_i}$ and $\Pi_i >0$ and $k \,>\,0$ are scalar gains. Each of $c$ and $k$ are defined to satisfy \eqref{eq:eq30} and \eqref{eq:eq31}
                            \begin{equation}
                             \label{eq:eq30}
                                k=\frac{c}{2\underline{\sigma}(Q)}
                            \end{equation}
                            \begin{equation}
                             \label{eq:eq31}
                                c\underline{\sigma}(Q)> \frac{1}{2}\phi_M\bar{\sigma}(P)\bar{\sigma}(A)
                            \end{equation}

           Therefore, the control node trajectory $x_0\left(t\right)$ is cooperative UUB and all nodes synchronize close to $x_0\left(t\right)$.
        \end{theorem}
        \textbf{Proof:}\\
               In view of \eqref{eq:eq28}, the error function in \eqref{eq:eq7} can be defined by
                \begin{equation}
                 \label{eq:eq32}
                    \dot{e} = \left(L+B\right)\left(\tilde{W}^{\top}\phi\left(x\right) + \alpha -c\epsilon + w -\underline{f}\left(x_0,t\right)\right)
                \end{equation}
                and from \eqref{eq:eq18}, the transformed error can be written as
                \begin{equation}
                 \label{eq:eq33}
                   \begin{aligned}
                    \dot{\epsilon} = & R\left(L+B\right)\left[\tilde{W}^{\top}\phi\left(x\right) + \alpha -c\epsilon + w -\underline{f}\left(x_0,t\right)\right] \\
                    &- R\,\,\dot{\varUpsilon}\,\varUpsilon^{-1}\,e
                    \end{aligned}
                \end{equation}
               Consider the following Lyapunov candidate function
                \begin{equation}
                 \label{eq:eq34}
                    V = \frac{1}{2}\epsilon^{\top}\,P\epsilon + \frac{1}{2}{ \rm Tr}\left\{\tilde{W}^{\top}F^{-1}\tilde{W}\right\}
                \end{equation}
               where $P$ is the positive and diagonal matrix defined in Lemma \ref{lemma1}. $F^{-1}$ is a block diagonal matrix defined in \eqref{eq:eq29} with $F :={\rm diag}\left\{F_i\right\}$. The derivative of \eqref{eq:eq34} is
                \begin{equation}
                 \label{eq:eq35}
                    \dot{V} = \frac{1}{2}\dot{\epsilon}^{\top}\,P\epsilon+\frac{1}{2}\epsilon^{\top}\,P\dot{\epsilon} + { \rm Tr}\left\{\tilde{W}^{\top}F^{-1}\dot{\tilde{W}}\right\}
                \end{equation}
               Let $P_1:=P\,R=R\,P$, such that $Q$ in \eqref{eq:eq266} is $Q=P_1\left(L+B\right)+\left(L+B\right)^{\top}P_1$. One can write
                \begin{equation}\label{eq:eq36}
                 \begin{aligned}
                    \dot{V} \leq & -\,c\frac{1}{2}\epsilon^{\top}\,Q\, \epsilon + \epsilon^{\top}\,P_1\,\left(L+B\right)(\alpha + w -\underline{f}\left(x_0,t\right))\\
                     &+ \epsilon^{\top}\,P_1\,\left(L+B\right)\tilde{W}^{\top}\phi\left(x\right) + { \rm Tr}\left\{\tilde{W}^{\top}F^{-1}\dot{\tilde{W}}\right\} \\
                     & - \epsilon^{\top}\,P_1\,\,\dot{\varUpsilon}\,\varUpsilon^{-1}\,e
                 \end{aligned}
                \end{equation}
        On the other hand $e= \varUpsilon\,S\left(\epsilon\right)$, which implies
        \begin{equation}
                 \label{eq:eq36b}
                 \begin{aligned}
                    \dot{V} = &-\frac{1}{2}\,c\epsilon^{\top}\,Q\,\epsilon + \epsilon^{\top}\,P_1\,\left(L+B\right)(\alpha + w -\underline{f}\left(x_0,t\right))\\
                     &+ \epsilon^{\top}\,P_1\,\left(L+B\right)\tilde{W}^{\top}\phi\left(x\right) + { \rm Tr}\left\{\tilde{W}^{\top}F^{-1}\dot{\tilde{W}}\right\} \\
                     &- \epsilon^{\top}\,P_1\,\,\dot{\varUpsilon}\,\,S\left(\epsilon\right)
                 \end{aligned}
                \end{equation}
        one should note that $\varLambda:=\varLambda\left(t\right)$ such that $\varLambda=-P_1\,\,\dot{\varUpsilon}$ is a positive definite diagonal matrix for $\forall t$ and $\lim\limits_{t \to \infty} \varLambda=0$
        \begin{equation}
                 \label{eq:eq36c}
                 \begin{aligned}
                    \dot{V} \leq &-\frac{1}{2}\,c\epsilon^{\top}\,Q\,\epsilon  + \epsilon^{\top}\,P\,\left(L+B\right)(\alpha + w -\underline{f}\left(x_0,t\right))\\
                     &+ \epsilon^{\top}\,P\,\left(L+B\right)\tilde{W}^{\top}\phi\left(x\right) + { \rm Tr}\left\{\tilde{W}^{\top}F^{-1}\dot{\tilde{W}}\right\}\\ &+ \epsilon^{\top}\,\varLambda\,\,\delta
                 \end{aligned}
                \end{equation}
       where $\bar{\delta}={ \rm max}\{\bar{\delta}_1, \,\ldots, \bar{\delta}_N\}$, equation \eqref{eq:eq36c} becomes
                \begin{equation}
                 \label{eq:eq38}
                 \begin{split}
                    \dot{V} \leq &-\frac{1}{2}\,c\epsilon^{\top}\,Q\,\epsilon + \epsilon^{\top}\,P_1\,\left(L+B\right)(\alpha + w -\underline{f}\left(x_0,t\right))\\
                     & +{\rm Tr}\left\{\tilde{W}^{\top}\,\left(F^{-1}\,\dot{\tilde{W}}+\phi\left(x\right)\,\epsilon^{\top}\,P_1\,\left(L+B\right)\right)\right\}\\
                     &+ \epsilon^{\top}\,\varLambda\,\,\bar{\delta} \,\underline{1}
                 \end{split}
                \end{equation}
        \begin{equation}
                 \label{eq:eq38b}
                 \begin{aligned}
                    \dot{V} \leq &-\frac{1}{2}\,c\epsilon^{\top}\,Q\,\epsilon + \epsilon^{\top}\,P_1\,\left(L+B\right)(\alpha + w -\underline{f}\left(x_0,t\right))\\
                     & +\, {\rm Tr}\left\{\tilde{W}^{\top}\,\left(F^{-1}\,\dot{\tilde{W}}+\phi\left(x\right)\,\epsilon^{\top}\,P_1\left(D+B\right)\right)\right\}\\
                     & -\,{\rm Tr}\left\{\tilde{W}^{\top}\,\phi\left(x\right)\,\epsilon^{\top}\,P_1\,A\right\}+ \epsilon^{\top}\,\varLambda\,\,\bar{\delta} \,\underline{1}
                 \end{aligned}
                \end{equation}
         Using the NN weight tuning law in \eqref{eq:eq29} and since $P$ and $D+B$ are diagonal, one has
          \begin{equation}
                   \label{eq:eq38c}
                   \begin{aligned}
                      \dot{V} \leq &-\frac{1}{2}\,c\epsilon^{\top}\,Q\,\epsilon + \epsilon^{\top}\,P_1\,\left(L+B\right)(\alpha + w -\underline{f}\left(x_0,t\right))\\
                       & + k\, {\rm Tr}\left\{\tilde{W}^{\top}\,\left(W-\tilde{W}\right)\right\}-k\,{\rm Tr}\left\{\tilde{W}^{\top}\,\phi\left(x\right)\,\epsilon^{\top}\,P_1\,A\right\}\\
                       &+ \epsilon^{\top}\,\varLambda\,\,\bar{\delta} \,{\bf \underline{1}}
                   \end{aligned}
                  \end{equation}

               Let $B_M = \alpha_M + w_M +F_M$, $z=\begin{bmatrix}
                                   \left\Vert \epsilon \right\Vert & \left\Vert \tilde{W} \right\Vert_F
                                   \end{bmatrix}^{\top}$\\$h=\begin{bmatrix}
                                   \bar{\sigma}\left(P_1\right) \left( \bar{\sigma}\left(L+B\right)B_M\right) + \bar{\delta}\bar{\sigma}\left(\varLambda\right)  & kW_M
                                   \end{bmatrix}^\top$\\
                $H=\begin{bmatrix}
                                      \frac{1}{2}\,c\,\underline{\sigma}(Q)  & -\frac{1}{2}\phi_M\bar{\sigma}\left(P_1\right)\bar{\sigma}\left(L+B\right) \\
                                      -\frac{1}{2}\phi_M\bar{\sigma}\left(P_1\right)\bar{\sigma}\left(L+B\right) & k
                                     \end{bmatrix}$
                
                hence, \eqref{eq:eq38c} can be expressed as follows
                \begin{equation}
                 \label{eq:eq41}
                 \begin{split}
                    \dot{V} \leq &-z^{\top}\,H\,z + h^{\top}\,z
                 \end{split}
                \end{equation}
                one can define $\dot{V} \leq 0$ if and only if $H$ is positive definite and
                \begin{equation}
                 \label{eq:eq42}
                 \begin{split}
                    \left\Vert z \right\Vert > & \frac{\left\Vert h \right\Vert}{\underline{\sigma}\left(H\right)}
                 \end{split}
                \end{equation}
               according to \eqref{eq:eq34}, we have
                \begin{equation}
                 \label{eq:eq43}
                 \begin{aligned}
                    \frac{1}{2}\underline{\sigma}\left(P_1\right)\left\Vert \epsilon \right\Vert^2 + \frac{1}{2\Pi_{\rm max}}\left\Vert \tilde{W} \right\Vert_F^2 \leq V 
                    \leq &\frac{1}{2}\bar{\sigma}\left(P_1\right)\left\Vert \epsilon \right\Vert^2 \\
                    &+ \frac{1}{2\Pi_{\rm min}}\left\Vert \tilde{W} \right\Vert_F^2
                 \end{aligned}
                \end{equation}
                or
                 \begin{equation}
                  \label{eq:eq44}
                  \begin{split}
                    & \frac{1}{2}\begin{bmatrix}
                    \left\Vert \epsilon \right\Vert & \left\Vert \tilde{W} \right\Vert_F
                    \end{bmatrix}
                    \begin{bmatrix}
                      \underline{\sigma}\left(P_1\right)  & 0\\
                      0 & \frac{1}{\Pi_{\rm max}}
                     \end{bmatrix}
                    \begin{bmatrix}
                    \left\Vert \epsilon \right\Vert \\
                    \left\Vert \tilde{W} \right\Vert_F
                    \end{bmatrix} \leq V \\
                    &\leq
                     \frac{1}{2}\begin{bmatrix}
                    \left\Vert \epsilon \right\Vert & \left\Vert \tilde{W} \right\Vert_F
                    \end{bmatrix}
                    \begin{bmatrix}
                      \bar{\sigma}\left(P_1\right)  & 0\\
                      0 & \frac{1}{\Pi_{\rm min}}
                     \end{bmatrix}
                    \begin{bmatrix}
                    \left\Vert \epsilon \right\Vert \\
                    \left\Vert \tilde{W} \right\Vert_F
                    \end{bmatrix}
                  \end{split}
                 \end{equation}
                with $\Pi_{\rm min}$ and $\Pi_{\rm max}$ are minimum and maximum values of $\Pi_{i}$, respectively. Define the appropriate variables in \eqref{eq:eq44} such that \eqref{eq:eq44} is written as
                 \begin{equation*}
                         \begin{split}
                            \frac{1}{2}z^{\top}\,\underline{S}\,z \leq V \leq \frac{1}{2}z^{\top}\,\bar{S}z
                         \end{split}
                        \end{equation*}
                which is equivalent to
                \begin{equation}
                 \label{eq:eq45}
                 \begin{split}
                    \frac{1}{2}\underline{\sigma}(\underline{S})\left\Vert z \right\Vert^2 \leq V \leq \frac{1}{2}\bar{\sigma}(\bar{S})\left\Vert z \right\Vert^2
                 \end{split}
                \end{equation}
                hence, one can find the following relation
                \begin{equation}
                 \label{eq:eq46}
                 \begin{split}
                    V > \frac{1}{2}\bar{\sigma}(\bar{S})\frac{\left\Vert h \right\Vert^2}{\underline{\sigma}^2\left(H\right)}
                 \end{split}
                \end{equation}

            selecting $k$ as in \eqref{eq:eq30} leads to
            {\begin{equation}
                     \label{eq:eq46b}
                     \begin{split}
                        {\underline{\sigma}\left(H\right)}= \frac{c\underline{\sigma}(Q) -\frac{1}{2} \phi_M\, \bar{\sigma}(P) \bar{\sigma}^2(A)}{2}
                     \end{split}
                    \end{equation}
           which is positive under condition\eqref{eq:eq29} and therefore $z$ is UUB. Consequently, the sufficient conditions for \eqref{eq:eq42} are

                \begin{equation}
                 \label{eq:eq49}
                 \begin{split}
                    \left\Vert \epsilon \right\Vert > \frac{\bar{\sigma}\left(P_1\right) \left( \bar{\sigma}\left(L+B\right)B_M\right) + \bar{\delta}\,\bar{\sigma}\left(\varLambda\right) + kW_M}{\underline{\sigma}\left(H\right)}
                 \end{split}
                \end{equation}
                \begin{equation}
                 \label{eq:eq50}
                 \begin{split}
                    \left\Vert \tilde{W} \right\Vert_F > \frac{\bar{\sigma}\left(P_1\right) \left( \bar{\sigma}\left(L+B\right)B_M\right) + \bar{\delta}\,\bar{\sigma}\left(\varLambda\right) + kW_M}{\underline{\sigma}\left(H\right)}
                 \end{split}
                \end{equation}
                also from \eqref{eq:eq45}, we have
                \begin{equation}
                 \label{eq:eq51}
                 \begin{split}
                    \left\Vert z \right\Vert \leq \sqrt{\frac{2V}{\underline{\sigma}(\underline{S})}}, \hspace{10pt}\left\Vert z \right\Vert \geq \sqrt{\frac{2V}{\bar{\sigma}(\bar{S})}}
                 \end{split}
                \end{equation}
                consequently, equation \eqref{eq:eq41} can be written as
                \begin{equation}
                 \label{eq:eq52}
                    \dot{V} \leq  -\beta_1\,V+ \beta_2\,\sqrt{V}
                \end{equation}
                with $\beta_1 := \frac{2\underline{\sigma}\left(H\right)}{\bar{\sigma}(\bar{S})}$ and $\beta_2 := \frac{\sqrt{2}\left\Vert h \right\Vert}{\sqrt{\underline{\sigma}(\underline{S})}}$ which leads to
               \begin{equation}
                 \label{eq:eq54}
                    \sqrt{V} \leq \sqrt{V(0)} + \frac{\beta_2}{\beta_1}
                \end{equation}
             using \eqref{eq:eq45}, one has
             \begin{equation}
             \begin{aligned}
             \left\Vert \epsilon \right\Vert \, \leq \, \left\Vert z \right\Vert \leq & \sqrt{\frac{\bar{\sigma}(\bar{S})}{\underline{\sigma}(\underline{S})}} \sqrt{||e\left(t_0\right)||^2+\left\Vert \tilde{W}\left(t_0\right) \right\Vert_F^2}\\
             &+\frac{\bar{\sigma}(\bar{S})}{\underline{\sigma}(\underline{S})}\frac{\left\Vert r \right\Vert}{\underline{\sigma}\left(H\right)}\equiv B_{0}
             \end{aligned}
             \end{equation}
             Therefore, $\epsilon$ is $\mathcal{L}_{\infty}$ and contained for all $t>t_0$ in the compact set $\Omega_0=\{\epsilon\left(t\right)| \left\Vert \epsilon\left(t\right) \right\Vert\leq B_0\}$.
             These results hold under the assumption that $x \in \Omega_x$ as per Property \ref{property1} for all $t \geq t_0$; and therefore, we need to establish that the proposed control law and the initial conditions do not force $x$ to get out of the compact set $\Omega_x$. As such, the analysis follows \cite{ioannou2003} and \cite{das_distributed_2010} with slight modification. Let $x_0 \in \Omega_{x_0}$. Suppose $\tilde{W} \in \Omega_W$ . Then (5) shows that
             \begin{equation*}
              \left\Vert x\left(t\right) \right\Vert \leq \frac{1}{\underline{\sigma}\left(L+B\right)}\left\Vert e\left(t\right) \right\Vert+\sqrt{N}\left\Vert x_0\left(t\right) \right\Vert
             \end{equation*}
             Owing to Property \ref{propertiesS},
             \begin{equation}
                           \left\Vert x\left(t\right) \right\Vert \leq \frac{1}{\underline{\sigma}\left(L+B\right)}\left\Vert S\left(r_{t_0}\right) \right\Vert+\sqrt{N}\left\Vert x_0\left(t\right) \right\Vert \equiv B_1
             \end{equation}
            The state is contained for all times $t \geq 0$ in a compact set $\Omega_{x} = \left\{x\left(t\right)| \left\Vert x\left(t\right) \right\Vert \leq B_1\right\}$  and this completes the proof. It should be remarked that $k$ and $\Pi_i$ are associated with nonlinearity compensation of adaptive estimate, $c$ controls the speed of convergence to the desired tracking output and has to be selected to satisfy \eqref{eq:eq30}.

     Finally, the algorithm of nonlinear single node dynamics such as equation \eqref{eq:eq1} can be summarized briefly as
            \begin{enumerate}
              \item Define the control design parameters such as $\bar{\delta}_i$, $\underline{\delta}_i$, $\rho_{i\infty}$, $\ell_{i}$, $\xi$, $\Pi_i$, $k$ and $c$.
              \item Evaluate local error synchronization from equation  \eqref{eq:eq4}.
              \item Evaluate the PPF from equation \eqref{eq:eq9}.
              \item Evaluate $r_i$ from equation \eqref{eq:eq17c}.
              \item Evaluate transformed error from equation  \eqref{eq:eq58b}.
              \item Evaluate control signal from equation \eqref{eq:eq27}.
              \item Evaluate the neuro-adaptive estimate from equation \eqref{eq:eq29}.
              \item Go to step 2.
           \end{enumerate}

                \begin{remark}
                    If we have $x_i \in \mathbb{R}^{n}, n > 1$, $u_i\in \mathbb{R}^{n}$, $f_i\left(x_i\right)\in \mathbb{R}^{n}$ and $w_i\in \mathbb{R}^{n}$,then the problem can be extended easily and the estimated weight will be written as
                 \begin{equation}
                  \label{eq:eq55}
                     \dot{\hat{W}}_i = -F_i\phi_i\epsilon_ir_i\left(p_i\left(d_i+b_i\right)\otimes \mathbb{I}_n\right) - kF_i\hat{W}_i
                 \end{equation}
                \end{remark}

\section{Simulation Results} \label{Sec5}
      {\bf Example 1:} Consider the problem in \cite{das_distributed_2010} with strongly connected network composed of 5 nodes and one leader connected to node 3 as in Fig. \ref{fig:fig2}.
      \begin{figure}[h!]
       \centering
       \includegraphics[scale=0.5]{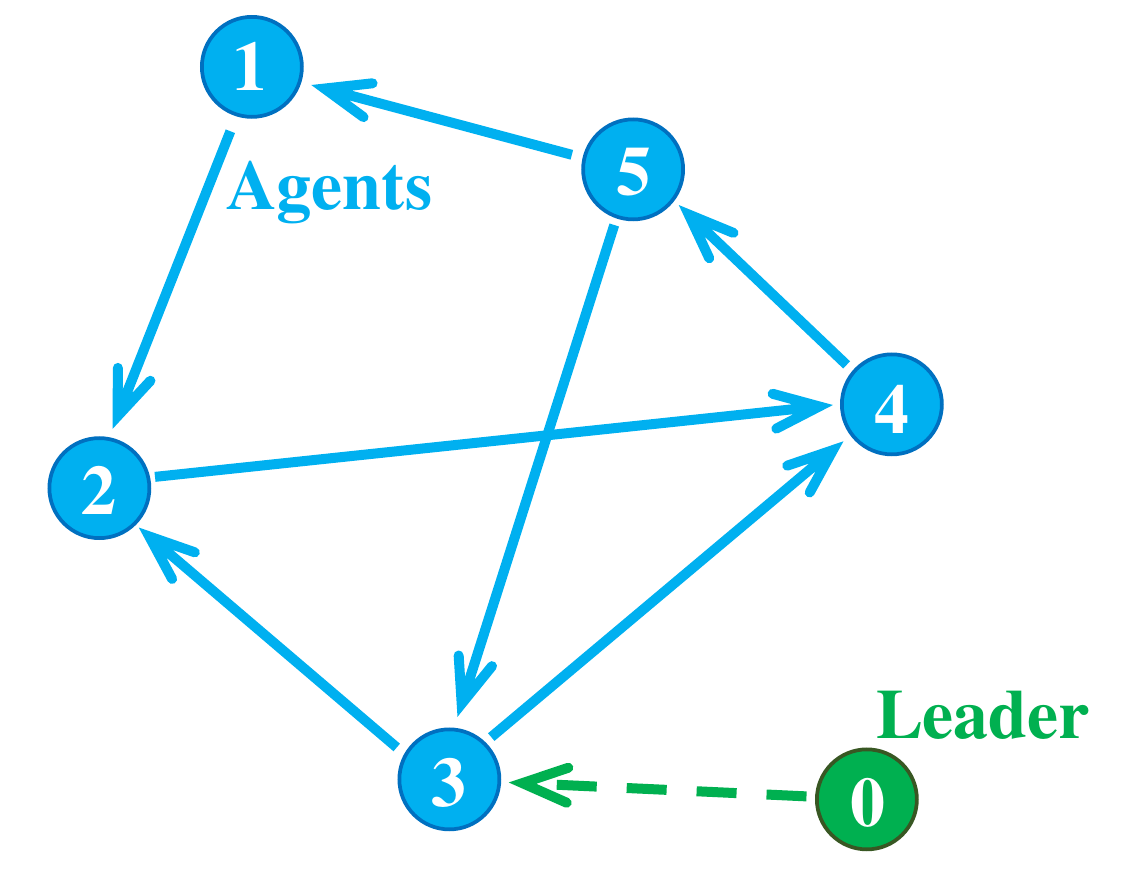}
       \caption{Strongly connected graph with one leader and five agents \cite{das_distributed_2010}.}
       \label{fig:fig2}
      \end{figure}

      The nonlinear dynamics of the graph are
       \begin{equation*}
          \begin{aligned}
            & \dot{x}_1 = x_1^3 + u_1 + a_1\left(t\right){\rm cos}\left(t\right)\\
            & \dot{x}_2 = x_2^2 + u_2 + a_2\left(t\right){\rm cos}\left(t\right)\\
            & \dot{x}_3 = x_3^4 + u_3 + a_3\left(t\right){\rm cos}\left(t\right)\\
            & \dot{x}_4 = x_4 + u_4 + a_4\left(t\right){\rm cos}\left(t\right)\\
            & \dot{x}_5 = x_5^5 + u_5 + a_5\left(t\right){\rm cos}\left(t\right)
          \end{aligned}
        \end{equation*}
       where $a_i\left(t\right)$, $i=1,\ldots,5$, represent random constants between 0 and 1 at each time instant and the leader dynamics is $\dot{x}_0 = f_0\left(x_0,t\right) = 0$ with desired consensus value = 2. The following parameters were used in the simulation
       $\rho_{\infty} = 0.05\times{\bf \underline{1}}_{1\times 5}$, $\rho_{0} = 7\times{\bf \underline{1}}_{1\times 5}$, $\ell = 7\times{\bf \underline{1}}_{1\times 5}$, $\Pi_i = 150 \, \mathbb{I}_{v_i\times v_i}$, $\bar{\delta} = 7\times{\bf \underline{1}}_{1\times 5}$, $\underline{\delta} = 1\times{\bf \underline{1}}_{1\times 5}$, $v_i = 3$, $c = 100$, $k = 0.8$, $\xi = 20$, $x_0 = 2{\rm cos}(0.8t)$ and
        $x_{i}(0) = \left[-2.5743 , -0.9814 , 1.2596 , 1.1472 , 2.5196\right]$.\\
       Fig. \ref{fig:fig3} shows the output performance, control signal and transformed error respectively for the proposed control algorithm using the transformed error in \eqref{eq:eq12b} instead of that in \eqref{eq:eq58b}. In all figures, $x_0$ denotes the output response of the leader node, $u_0$ denotes the control signal of the leader node and $x_i$ and $u_i$ are the output performance and control signal of the associated agent, respectively, for all $i=1,\ldots,5$. Fig. \ref{fig:fig3} illustrates the quality of tracking and how it falls between the constraints of the PPF. Also, it can be noticed that the error obeys prescribed error boundaries as it started within large set and ended within small predefined set. However, the control input and transformed error are completely oscillatory.

      \begin{figure}[h!]
       \centering
       \includegraphics[scale=0.27]{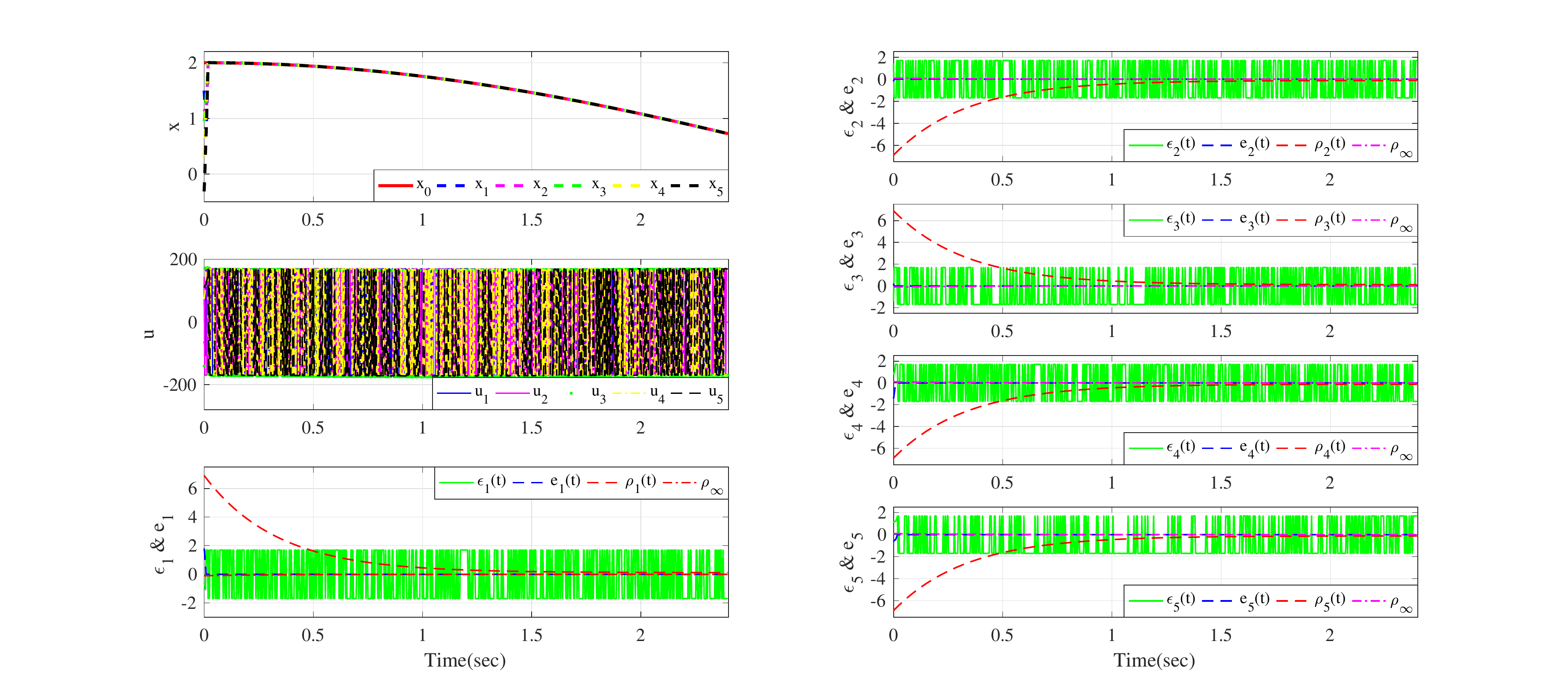}
       \caption{The output performance and control signal of example 1 using the transformed error in equation \eqref{eq:eq12b}.}
       \label{fig:fig3}
      \end{figure}

       Now, the proposed neuro-adaptive distributed control with PPF has been improved to handle the problem of oscillation in  control input and transformed error of Fig. \ref{fig:fig3} using the transformed error in \eqref{eq:eq58b}. Fig. \ref{fig:fig5} and \ref{fig:fig6} present the output performance, control signal and transformed error respectively using the proposed control algorithm with the smooth error function, $erf()$, in the transformed error \eqref{eq:eq58b}. The advantage of introducing such smooth function is reflected in the good tracking
       performance and with no oscillation in the behavior of transformed error and control signal response as illustrated on Fig. \ref{fig:fig6} and \ref{fig:fig5} respectively.
       Also, Fig. \ref{fig:fig5} compares the proposed neuro-adaptive distributed control with PPF and neuro-adaptive distributed control without PPF \cite{das_distributed_2010}. With same design parameters, initial conditions and leader response  $x_0 = 2{\rm cos}(0.8t)$, the proposed algorithm illustrates high tracking performance with low control effort, whereas the control in \cite{das_distributed_2010} shows instability.
      
      \begin{figure}[h!]
       \centering
       \includegraphics[scale=0.27]{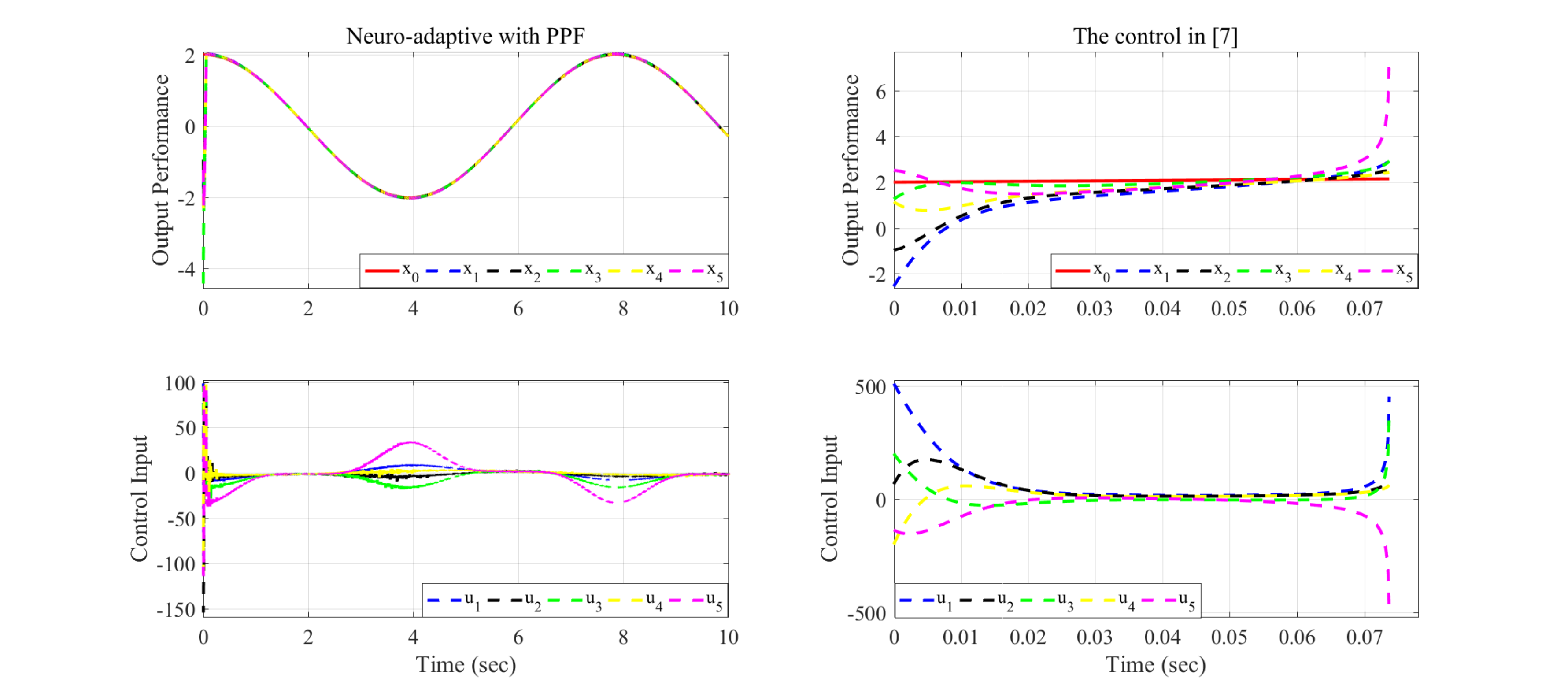}
       \caption{The output performance and control signal of Example 1 of the neuro-adaptive control with PPF using the transformed error in equation \eqref{eq:eq58b} compared to the neuro-adaptive control in \cite{das_distributed_2010}.}
       \label{fig:fig5}
      \end{figure}      

      \begin{figure}[h!]
       \centering
       \includegraphics[scale=0.27]{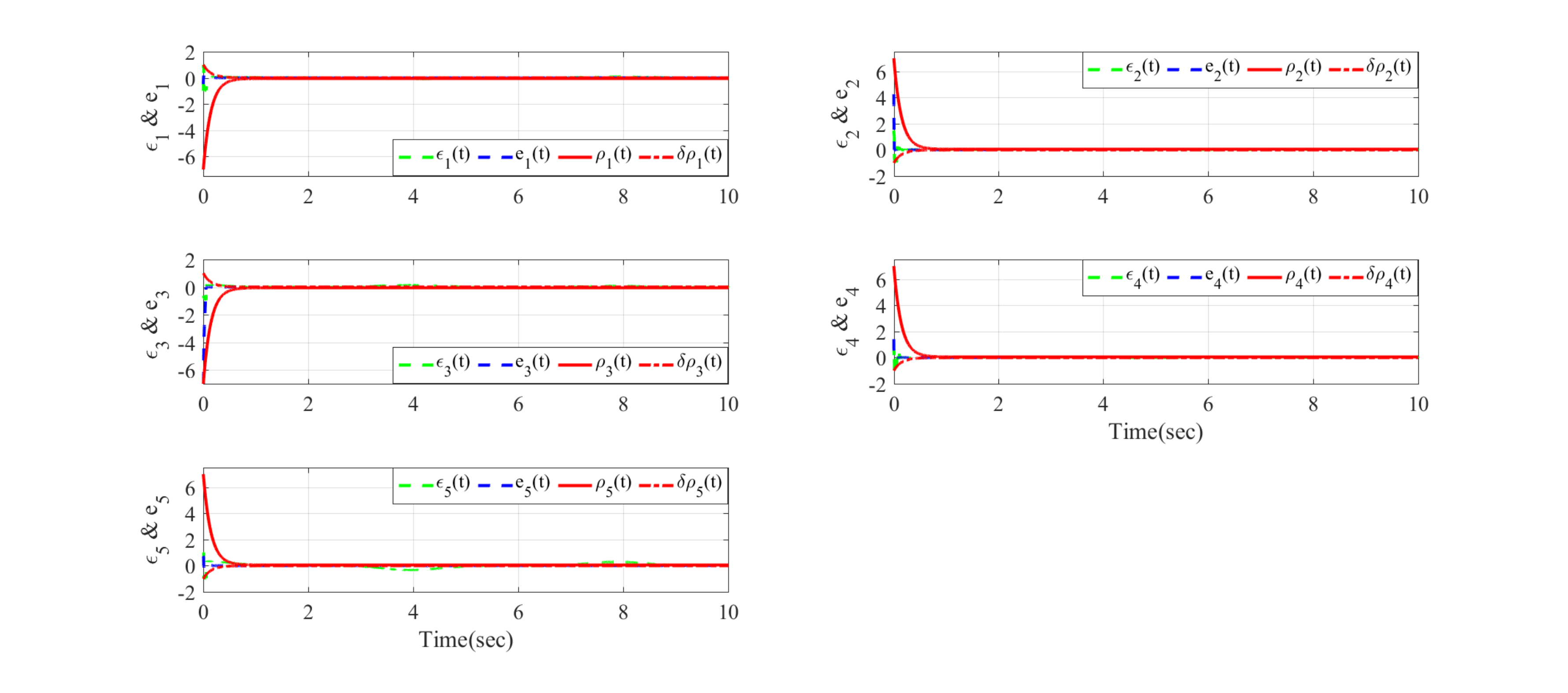}
       \caption{Error and transformed error of neuro-adaptive control with PPF of Example 1 using the transformed error in equation \eqref{eq:eq58b}.}
       \label{fig:fig6}
      \end{figure}

      {\bf Example 2:} Consider the same problem in Fig. \ref{fig:fig2} with 3 inputs and 3 outputs nonlinear systems. The nonlinear dynamics are defined by
       \begin{equation*}
          \begin{aligned}
             \dot{x}_{:,j} & = A\,x_{:,j} + B\,u_{:,j} + \theta_j\left(t\right) x_{:,j} + f_j\left(x_{:,j}\right) + D_j\left(t\right) \\
             y_{:,j} & = C\,x_{:,j}
          \end{aligned}
        \end{equation*}
        where
       \begin{equation*}
          \begin{aligned}
            A &= \begin{bmatrix}
            -20 & 22 & 0\\ 0 & 15 & 0 \\  0 &  0  & -3
            \end{bmatrix},
            B = \mathbb{I}_{3 \times 3},C = \mathbb{I}_{3 \times 3},\\
            f_j(x_{:,j}) &= \begin{bmatrix}
              a_{1,j}x_{3,j}x_{1,j}+0.2{\rm sin}\left(x_{1,j}a_{1,j}\right)\\ -a_{2,j}x_{1,j}x_{3,j}-0.2a_{2,j}{\rm cos}\left(a_{2,j}x_{3,j}t\right)x_{1,j}\\     a_{3,j}x_{1,j}x_{2,j}
            \end{bmatrix},\\
            D_j\left(t\right) &= \begin{bmatrix}
              1+b_{1,j}{\rm sin}\left(b_{1,j}t\right)\\
               1.2{\rm cos}\left(b_{2,j}t\right)\\
              {\rm sin}\left(0.5b_{3,j}t\right) + {\rm cos}\left(b_{3,j}t\right) - 1
            \end{bmatrix},
          \end{aligned}
        \end{equation*}
        {\scriptsize
       \begin{equation*}
            \theta_j^1\left(t\right) = \begin{bmatrix}
              3c_{1,j}{\rm sin}\left(0.5t\right)\\
              0.9{\rm sin}\left(0.2c_{2,j}t\right)\\
              0.5{\rm sin}\left(0.13c_{3,j}t\right)
            \end{bmatrix},\hspace{10pt}
                      \theta_j^2\left(t\right) = \begin{bmatrix}
                         2c_{1,j}{\rm sin}\left(0.4c_{1,j}t\right){\rm cos}\left(0.3t\right)\\
                        2.5{\rm sin}\left(0.3c_{2,j}t\right)+0.3{\rm cos}\left(t\right)\\
                         0.6c_{3,j}{\rm cos}\left(0.15t\right)
                      \end{bmatrix},
        \end{equation*}
      \begin{equation*}
                      \theta_j^3\left(t\right) = \begin{bmatrix}
                        0.7{\rm sin}\left(0.2c_{1,j}t\right)\\
                        1.0{\rm sin}\left(0.1c_{2,j}t\right)\\
                        1.5{\rm cos}\left(0.7c_{3,j}t\right)+1.6c_{3,j}{\rm sin}\left(0.3t\right)
                      \end{bmatrix},
        \end{equation*}
        }
               \begin{equation*}
                  \begin{aligned}
                    \theta_j\left(t\right) = \begin{bmatrix}\theta_j^1\left(t\right)&\theta_j^2\left(t\right)&\theta_j^3\left(t\right) \end{bmatrix}
                  \end{aligned}
                \end{equation*} 
      and
       \begin{equation*}
          \begin{aligned}
            a = \begin{bmatrix}
              1.5 & 0.5 & 0.7 & 1.3 & 0.7\\
              0.5 & 1.4 & 0.1 & 1.3 & 2.4\\
              2.8 & 1.4 & 0.6 & 0.7 & 0.6
            \end{bmatrix},\\
            b = \begin{bmatrix}
              0.5 & 1.5 & 1.1 & 1.6 & 0.3\\
              0.7 & 1.2 & 1.3 & 0.5 & 0.3\\
              1.1 & 1.4 & 1.6 & 0.6 & 1.0
            \end{bmatrix},\\
            c = \begin{bmatrix}
              1.5 & 2.5 & 0.5 & 1.7 & 0.7\\
              0.5 & 1.7 & 1.1 & 0.3 & 0.4\\
              0.8 & 0.4 & 2.2 & 0.9 & 1.4
            \end{bmatrix},
          \end{aligned}
        \end{equation*}
      The leader dynamics is $\dot{x}_0 = f_0\left(x_0,t\right) = \left[0,0,0\right]^{\top}$ with initial conditions $x_0(0)=\left[1.5,2.7,3.5\right]^{\top}$ equal to desired consensus and number of neurons $v_i = 6$. The control parameters of the problem were defined as $\rho_{\infty} = 0.05\times{\bf \underline{1}_{3\times 5}}$, $\rho_{0} = 7\times{\bf \underline{1}_{3\times 5}}$, $\ell = 7\times{\bf \underline{1}_{3\times 5}}$, $\Pi_i = 150\mathbb{I}_{v_i\times v_i}$, $\bar{\delta} = 7\times{\bf \underline{1}_{3\times 5}}$, $\underline{\delta} = 1\times{\bf \underline{1}_{3\times 5}}$, $c = 100$, $k = 0.8$, $\xi = 50$. Initial conditions of
             \begin{equation*}
                \begin{aligned}
                 x_1(0) &= [0.8678, -0.5058,1.5501]^{\top}\\
                 x_2(0) &=[-0.7777,0.3450, -2.5642]^{\top}\\
                 x_3(0) &=[2.1855, -2.5301, -1.3017]^{\top}\\
                 x_4(0) &=[1.4068, -3.3671, 0.8756]^{\top}\\
                 x_5(0)&=[1.0541, 0.2989, -1.4818]^{\top}
                \end{aligned}
             \end{equation*}
      The robustness of the proposed controller against time varying uncertain parameters, time varying external disturbances, and high nonlinearities are tested in this example. Fig. \ref{fig:fig7} shows the output performance and the control input of the proposed controller for this MIMO case. The tracking errors and their associated transforms are depicted in Fig. \ref{fig:fig9}, \ref{fig:fig10} and \ref{fig:fig11}. The control signal is smooth and the dynamic performance satisfy the prescribed dynamic constraints.

      \begin{figure}[h!]
       \centering
       \includegraphics[scale=0.27]{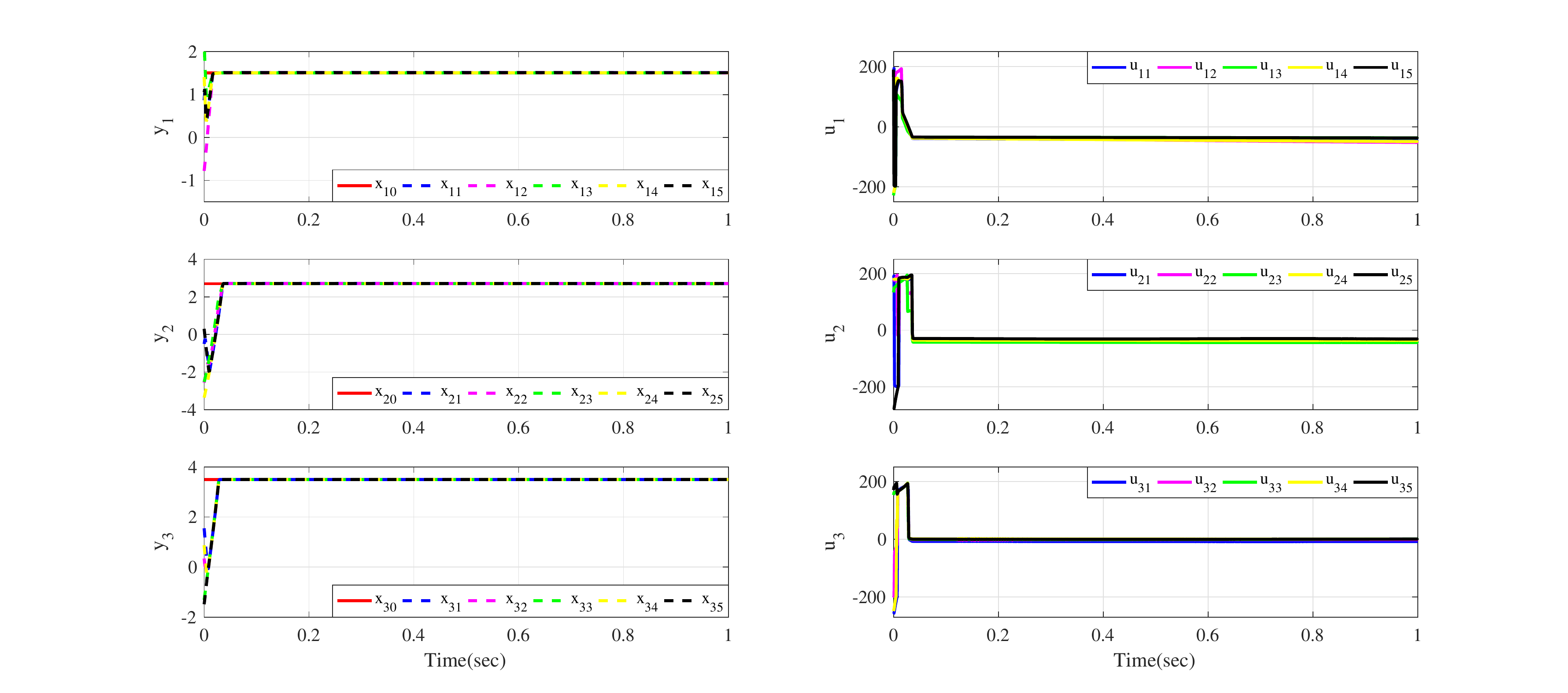}
       \caption{Neuro-adaptive control with PPF output performance of Example 2 using the transformed error in equation \eqref{eq:eq58b}.}
       \label{fig:fig7}
      \end{figure}

      \begin{figure}[h!]
       \centering
       \includegraphics[scale=0.27]{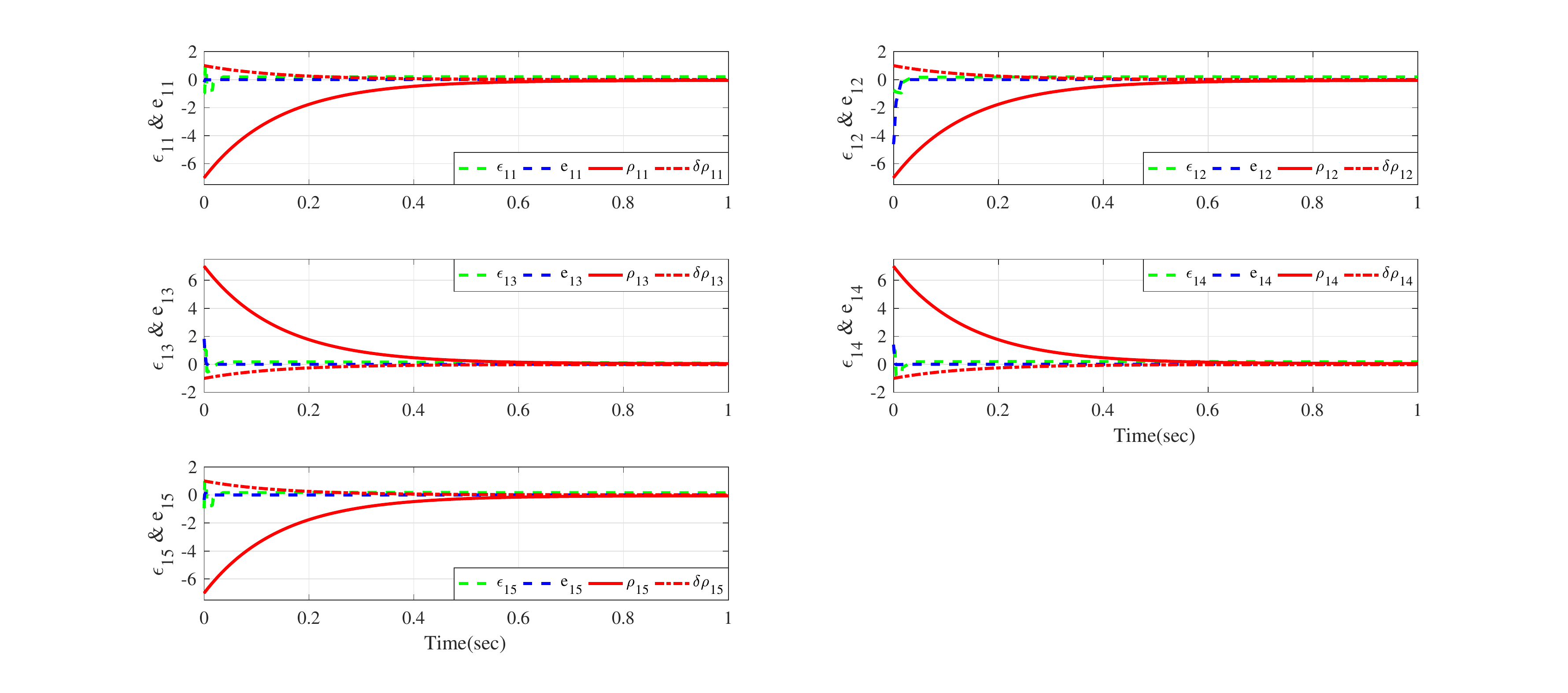}
       \caption{Error and transformed error of output 1 of example 2 using the transformed error in equation \eqref{eq:eq58b}.}
       \label{fig:fig9}
      \end{figure}

      \begin{figure}[h!]
       \centering
       \includegraphics[scale=0.27]{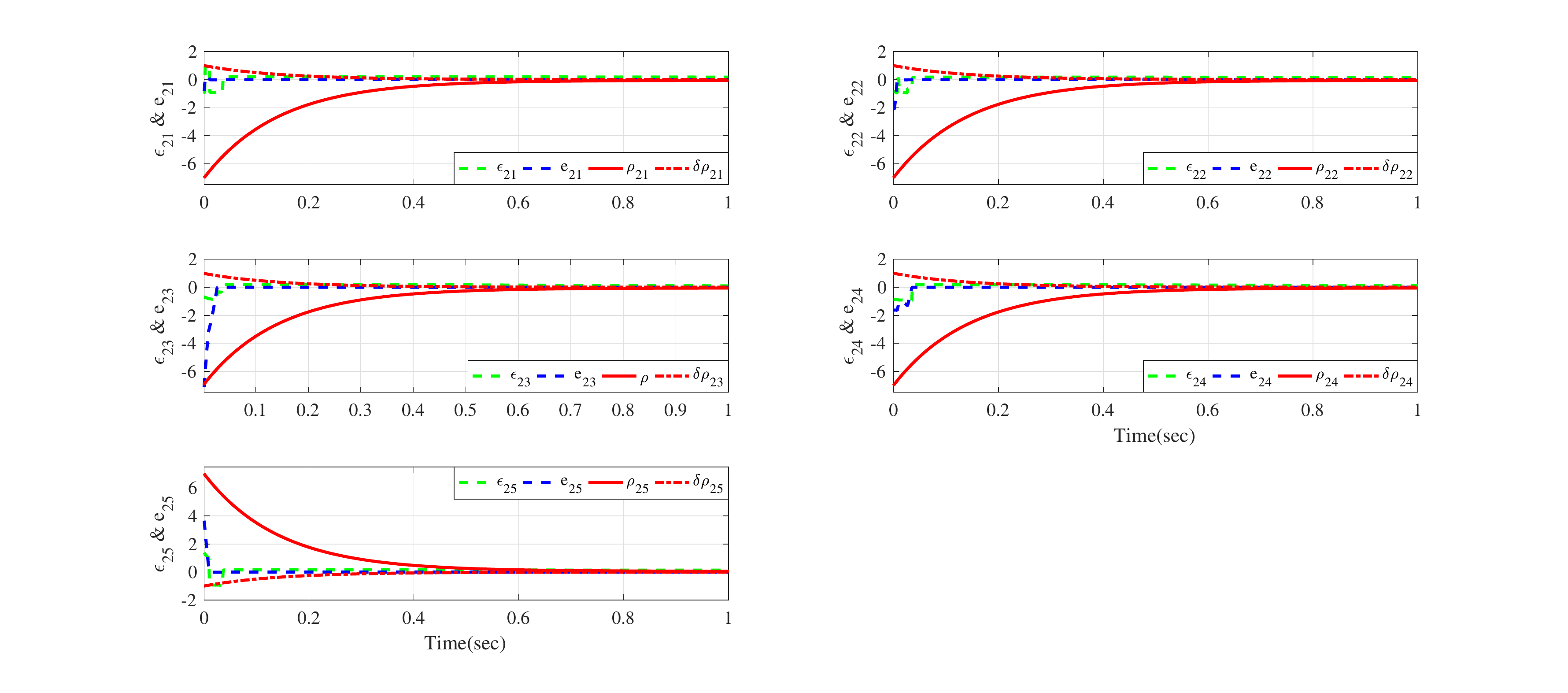}
       \caption{Error and transformed error of output 2 of Example 2 using the transformed error in equation \eqref{eq:eq58b}.}
       \label{fig:fig10}
      \end{figure}

      \begin{figure}[h!]
       \centering
       \includegraphics[scale=0.27]{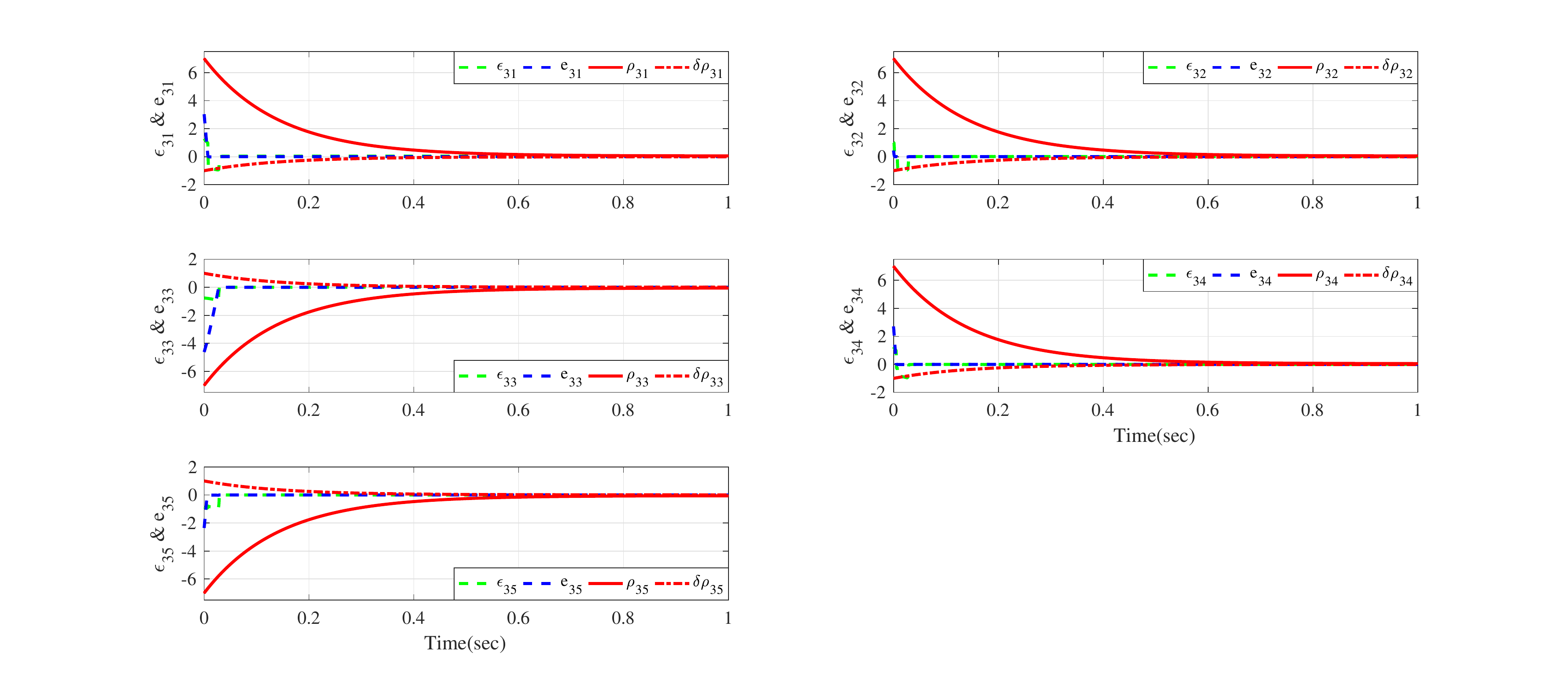}
       \caption{Error and transformed error of output 3 of Example 2 using the transformed error in equation \eqref{eq:eq58b}.}
       \label{fig:fig11}
      \end{figure}

\section{Conclusion}  \label{Sec6}
    A distributed cooperative neuro-adaptive control has been proposed based on prescribed performance for multi-agent systems. Neural network has been used to estimate system nonlinearities. The controller is developed for strongly connected network. The control signal is chosen properly to ensure  graph stability as well as track the leader trajectory with an error subject to dynamic prescribed constraints. The controller successfully enables the agents to synchronize to the desired trajectory in a very short time with small residual error. Simulation examples have been presented  by considering highly nonlinear heterogeneous systems with time varying parameters, uncertainties and disturbances. Only first order nonlinear systems have been considered. Higher order nonlinear systems will be considered in the future.

\section*{Acknowledgment}                               
The author(s) would like to acknowledge the support provided by
King Abdulaziz City for Science and Technology (KACST) through the
Science and Technology Unit at King Fahd University of Petroleum
and Minerals (KFUPM) for funding this work through project No.
09-SPA783-04  as part of the National Science, Technology and
Innovation Plan.  

\ifCLASSOPTIONcaptionsoff \newpage{}\fi



\bibliographystyle{ieeetr}        
\bibliography{bib_MCC_PPF}           

%

%
%
%




\end{document}